\documentclass{article}
\usepackage[a4paper, total={6in, 8in}]{geometry}
\usepackage{graphicx}
\usepackage{float}
\usepackage{natbib}
\usepackage{subcaption}
\usepackage{booktabs}
\usepackage{adjustbox}
\usepackage{multirow}
\usepackage{xcolor}
\usepackage[normalem]{ulem}
\usepackage{lipsum}
\usepackage{tablefootnote}

\newcommand{\companyModel}{Company Model}
\newcommand{\CM}{CM}

\newcommand{\histShift}{Historical Shifted}
\newcommand{\HS}{HS}

\newcommand{\myparagraph}[1]{\paragraph{#1.}}

\makeatletter
\newcommand\blfootnote[1]{%
  \begingroup
  \renewcommand{\@makefntext}[1]{\noindent\makebox[1.8em][r]#1}
  \renewcommand\thefootnote{}\footnote{#1}%
  \addtocounter{footnote}{-1}%
  \endgroup
}
\makeatother

\title{The impact of machine learning forecasting on strategic decision-making for Bike Sharing Systems}
\author{\it E. Angelelli$^1$, A. Mor$^{2*}$, C. Orsenigo$^2$, M.G. Speranza$^1$, C. Vercellis$^2$\blfootnote{ORCID. Enrico Angelelli: 0000-0003-1847-9331. Andrea Mor: 0000-0001-6131-7229. Carlotta Orsenigo: 0000-0001-8688-414X. M. Grazia Speranza: 0000-0002-8893-5227. Carlo Vercellis: 0000-0002-5020-3688.} \\
{\small \it $^1$Department of Economics and Management}\\
{\small \it University of Brescia, Italy}\\
{\small \{enrico.angelelli, grazia.speranza\}@unibs.it}\\
{\small \it $^2$Department of Management, Economics and Industrial Engineering}\\
{\small \it Politecnico di Milano, Italy}\\
{\small \{andrea.mor, carlotta.orsenigo, carlo.vercellis\}@polimi.it}\\
{\small $^*$Corresponding author}\\
}
\date{}

\begin{document}

\maketitle


\section*{Abstract}

In this paper, machine learning techniques are used to forecast the difference between bike returns and withdrawals at each station of a bike sharing system. The forecasts are integrated into a simulation framework that is used to support long-term decisions and model the daily dynamics, including the relocation of bikes. 
We assess the quality of the machine learning-based forecasts in two ways. Firstly, we compare the forecasts with alternative prediction methods. Secondly, we analyze the impact of the forecasts on the quality of the output of the simulation framework. The evaluation is based on real-world data of the bike sharing system currently operating in Brescia, Italy.

\section*{Practitioner summary}

Effective strategic planning is fundamental for the long-term sustainability of Bike Sharing Systems. For operators and policy makers, a primary challenge lies in making informed decisions about infrastructure, fleet sizing, and resource allocation. This study explores how advanced Machine Learning techniques can support these high-level decisions by providing more accurate demand forecasting compared to traditional methods.

Leveraging extensive real-world data from the Bicimia system in Italy, we compared various forecasting models over a long-term horizon. This time-frame allows the study to account for significant daily and seasonal fluctuations in user demand. Our analysis suggests that a ``local'' approach, tailoring specific ML models to the unique usage pattern of each individual station, tends to provide better accuracy than ``global'' system-wide models. We integrated these forecasts into a simulation framework that explicitly considers the dynamic relocation of bikes. This allowed us to assess how improved forecasts impact service quality across different scenarios, including examples focused on fleet allocation with a varying number of relocation vehicles.

The results provide useful insights for practitioners: the proposed approach was associated with a 5.4\% reduction in unsatisfied demand compared to current practices, potentially enhancing service reliability for the public. Beyond performance metrics, the study indicates that utilizing advanced predictive models may lead to more efficient decisions regarding fleet management. These findings suggest that ML can serve as a valuable strategic asset, helping authorities transition toward more data-driven planning to improve the overall resilience and accessibility of urban transportation services.

\section{Introduction}

As the population tends to concentrate in urban areas and acting against climate change becomes more compelling, the availability of sustainable means of transportation becomes more and more urgent.
Bike Sharing Systems (BSSs) represent a way for public administrations (and private companies) to offer an environmentally friendly and  capillary travel method. Beside being a standalone option, BSSs are also beneficial when considered in synergy with other means of transportation (typically of higher capacity, such as light trains lines), covering the first and last mile of travel demand. 

Two main paradigms have emerged over time regarding the dynamics of bike  withdrawals and returns: dock-based and free-floating BSSs, with an in-between solution based on dockless hubs/bays systems. In the former, stations are scattered in the operating area and provide docks for bikes to be locked to or unlocked from. This has the downside of requiring specialized hardware for the bike stations, which increases the installation cost, while also constraining the locations where a bike  can be withdrawn or returned. Upsides include the low cost of the bikes, which typically do not require specific hardware for remote locking and tracking by the company, and the potentially low relocation cost, as relocation only takes place among the stations.
In the free-floating BSSs, bikes are scattered in the operating area and  users are able to start their trip in any location where a bike is available and end it close to their destination. This paradigm removes the need for specialized hardware for the stations but requires bikes to be customized in order to be traceable all time, for both users and operators needs. Moreover, it can increase the relocation cost, as relocating vehicles have to collect bikes potentially scattered in any point of the operating area (as opposed to the finite number of stations in a station-based BSS).
A mid-way option is based on dockless hubs/bays systems, where traceable bikes must be picked up and returned to stations which are geographically defined areas, e.g., a rectangular area in a city square. The presence of these areas allows for an easy relocation of bikes while the traceability of the bikes allows these systems to forego hardware costs of dock-based stations.
The latest edition of the Meddin Bike-Sharing World Map report (see \cite{meddin2024}) reports that the number of systems relying on stations (i.e., dock-based systems and dockless hubs/bays systems) is about 74\% of the total.
Henceforth, we will focus on these systems, which we will call ``station-based'' systems.

When designing and operating a BSS, several decisions spanning the strategic, tactical, and operational levels must be made and reassessed over time \citep{shui2020review,caggiani2021toward,zhang2023relocation}. At the strategic level, these decisions include the station network design (i.e., the stations location and capacity), the bike fleet sizing, 
 and the relocation strategy.
At the tactical and operational levels,
operating a BSS requires to deal, among others, with the imbalance of transportation demand both in time and space. This issue is typically addressed in two ways or a combination of them. The first is the active relocation of bikes through vehicles, either during the night when few to no people use the system (see \cite{raviv2013static})
or during the day considering the ongoing and forecasted use
(see \cite{ghosh2017dynamic}). The second is the management and shaping of demand, typically performed through systems of incentives such as varying prices (see \cite{haider2018inventory}).
Simulation tools play a pivotal role in the planning and operating of BSSs (see \cite{angelelli2022simulation}). By allowing decision-makers to model the movement of bikes between (existing or planned) stations, user demand and relocation strategies, these tools enable a comprehensive evaluation of strategic, tactical, and operational decisions, providing insights through performance metrics such as service quality, system efficiency, and operational costs. To effectively simulate a BSS,  users demand, in terms of bike withdrawals and returns, needs to be carefully forecasted. Such forecasts are used to plan the relocation of bikes.

The scope of this paper is to study the impact of methods based on Machine Learning (ML), used to forecast users demand, on the simulation of a station-based BSS aimed at supporting  strategic and tactical decision-making. 
To this end, three well-established ML methods are chosen to forecast the half-hourly net demand of bikes, that is the difference between returns and withdrawals, for each of the stations of a station-based BSS. 
The one-year look-ahead horizon considered in this paper far exceeds those presented in the literature, especially when compared to the half-hour time unit of the forecast. 
This long look-ahead period of time is required by the strategic and tactical nature of the decisions we aim to address, while the half-hour time unit allows for the demand forecast to guide operational decisions regarding bike relocation.
Two approaches are presented and tested to generate the forecasting model. In the first, a unique model is trained for all stations (``global'' approach), while in the second one model is trained for each station (``local'' approach).
To study the impact of the forecasting methods on long-term decision making, the forecasts are integrated into a simulation framework adapted from \cite{angelelli2022simulation}. Such framework models daily system dynamics and uses the predicted demand to guide dynamic bike relocation throughout the day.
The impact of ML-based forecasting on decision-making is assessed on the real-world case of the Bicimia BSS, operating in Brescia, Italy. We compare the ML-based forecasting methods with the current approach adopted by the company, and show how different forecasting methods lead to different strategic and tactical decisions (see \cite{mor2026optimization} for a discussion on the use of machine learning models and classical tools for decision support).
An example is presented for the allocation of relocation vehicles to work shifts, demonstrating how the choice of forecasting model influences scheduling outcomes and service quality.

The rest of this paper is organized as follows. 
In Section \ref{sec:litrev} a review of the literature on 
forecasting alone and forecasting for decision-making in BSSs is presented. 
In Section \ref{sec:forecast} the explanatory features used to forecast the user behavior in the BSS are reported, together with the ML algorithms adopted for the same task. The two approaches considered to build the forecasting models
are also defined.
In Section \ref{sec:simfw} the simulation framework used in this work is illustrated. 
In Section \ref{sec:bsmob} the real-world case of Bicimia is described. 
In Section \ref{sec:results} the computational results are presented.  In particular, we first assess the quality of the forecasting based on ML methods. Then, we evaluate the impact of the forecasting methods on the output of the simulation.  
Finally, conclusions are drawn in Section \ref{sec:conclusions}.

\section{Literature review}\label{sec:litrev}

In this section, the literature related to this paper is reviewed. 
First, we analyze how the issue of forecasting the demand in station-based BSSs has been addressed. 
Then, existing studies focused on the forecasting problem combined with that of decision-making for BSSs are reviewed.

\myparagraph{Forecasting for BSSs}

Forecasting for BSSs has been tackled in different ways depending on a variety of factors. These include the forecasting area (e.g., the entire system, clusters of stations, or single stations), the target of the forecast (e.g., total number of trips, number of withdrawals and returns), the forecast time unit, 
the forecast look-ahead horizon (i.e., the number of time units in the future for which a forecast is provided), and the forecasting method.
A recent review of the literature on machine learning methods in shared mobility is presented in \cite{teusch2023systematic}, while one on deep learning approaches to forecast BSS usage is presented in \cite{jiang2022bike}.

In the following, we report some of the most relevant contributions. A summary of these contributions is given in Table \ref{tab:ref}. 

As mentioned, a factor defining the forecasting approach is its spatial granularity, that is, the definition of the area being used to aggregate user demand. Depending on the application, and sometimes on the quality and sparsity of the data, different definitions of this dimension have been adopted. 
The entire area in which the system operates is often chosen when the aim is to extract the elements that influence the use of the system (e.g., weather-related, see \cite{ermagun2018urban}) or when other factors are considered, such as model complexity (i.e., the number of parameters of the model to be calibrated, see \cite{cantelmo2020low}). 
Partitions of the BSS operating area are, instead, often used to discretize the withdrawal and return area in free-floating BSS or to aggregate the data to deal with data sparsity (see, for instance, \cite{li2021short}).
Forecasting at the station level is the most appropriate one when decisions have to be made or revised for a station-based BSS. 
This level of granularity enables more actionable insights for operational planning, such as dynamic relocation of bikes, as it allows  the capturing of localized usage patterns.

Depending on the area chosen to aggregate demand, different quantities can be defined as the target of the forecasting task.
When system-wide forecasting is considered, the target is typically the total usage demand, that is, the number of bike trips in the entire system within a time unit.
When individual stations, or partitions of the operating area, are considered instead, two types of target are typically analyzed. The first type considers each station/partition individually. Possible targets include the number of bikes withdrawn and returned, or a function of these quantities (e.g., the difference between withdrawals and returns). An alternative target of this type is the stock of bikes, that is, the number of bikes present at the station in a certain time unit. 
A second type of target considers the number of trips  between each pair of stations/partitions. This latter case is reported in Table \ref{tab:ref} with the ``(OD)'' addendum to the ``Area'' column, indicating that the target quantity is referred to an Origin-Destination (OD) pair. 

Different time units have been considered in the literature when tackling the forecasting of the demand in a BSS, 
with the most common value being 1 hour.
As can be seen from Table \ref{tab:ref}, with respect to the look-ahead horizon, the most common approach is to forecast only the next time unit. When additional time periods are considered, the number of look-ahead time units 
varies from 12 \citep{yoon2012cityride} to 30 \citep{kaltenbrunner2010urban}.
While these cases already pose considerable challenges, they remain primarily oriented towards short-term operational decisions. They do not fulfill the requirements of forecasts that feed a simulation framework aimed at supporting long-term tactical and strategic decision-making, where longer time horizons are considered.

The adopted methods also vary widely, encompassing statistics-derived approaches such as Moving Average, Autoregressive Models, and Structural Equation Models, as well as Machine Learning–based techniques, including ensemble methods like Random Forest and Gradient Boosting, and Deep Learning–based approaches, such as Convolutional and Recurrent Neural Networks.

\begin{table}[H]
	\centering
	\begin{adjustbox}{max size={1\textwidth}{1\textheight}}
	\begin{tabular}{l|ccccp{7cm}}
		Reference & Area & Target & Time unit & Look-ahead (if present) & Method(s) \\ \hline
		\cite{liu2019multi} & Station & Stock & 1 Minute & 20 Min. & Long Short-Term Memory (LSTM)\\
		\cite{kaltenbrunner2010urban} & Station & Stock & 2 Min. & 1 Hours & Auto-Regressive Moving Average, other ad-hoc methods \\ 		
        \cite{yoon2012cityride} & Station & Stock & 5 Min. & 5,60 Min. & Modified Auto-Regressive Integrated Moving Average \\
		\cite{wang2018short} & Station & Stock & 5,10 Min. & - & LSTM, Gated Recurrent Unit (GRU), Random Forest (RF) \\
		\cite{almannaa2020dynamic} & Station & Stock & 15,30,45,60,120 Min. & - & Dynamic Linear Models \\
		\cite{collini2021deep} & Station & Stock & 15,30,45,60 Min. & - & Bidirectional LSTM \\	 
		\cite{kim2019graph} & Station & Withd. & 1 Hour & - & Graph Convolutional Neural Network \\
		\cite{rixey2013station} & Station & Withd. & 1 Month & - & Multivariate Linear Regression \\
		\cite{yang2016mobility} & Station & Withd. and return & 2,5,10 Min. & 30 Min. & RF, Historical Average, ARMA, other ad-hoc methods \\ 		
        \cite{luo2021predicting} & Station & Withd. and return & 5 Min. & 5,10,15,20 Min. & Local Spectral Graph Convolution-LSTM \\
		\cite{ma2022short} & Station & Withd. and return & 15,30,45,60 Min. & - & Spatio-Temporal Graph Attention-LSTM \\ 
        \cite{mrazovic2018deep} & Station & Withd. and return & 30 Min. & 4 Hours & LSTM \\ 	
		\cite{chen2020predicting} & Station & Withd. and return & 1 Hour & - & Shared Forward NN+Bidirectional RNN+Soft Attention Mechanism \\
		\cite{sohrabi2020real} & Station & Withd. and return & 1 Hour & - & Generalized Extreme Value Count Models \\ 
		\cite{boufidis2020development} & Station & Withd. and return & 1,2,3 Hours & - & Gradient Boosting (GB), XGBoost, RF, NN \\
		\cite{zhou2018markov} & Station & Withd. and return & 1 Day & - & Markov chains\\
        \cite{torres2024forecasting} & Station & Withd. and return & 1 and 3 Hours & - & RF, GB, NN \\ 
		\cite{sohrabi2021dynamic} & Station & (Net of) Withd. and return & 15 Min. & 15 Min. to 4 Hours & Pattern recognition + $k$-Nearest Neighbor\\ 
		\cite{zhou2019reliable} & Station (OD) & Net demand & 1 Hour & - & STW+M \\
		\cite{ranaiefar2016bike} & Station (OD) & Travels, Withd. and return & 1 Month & - & Structural Equation Modeling \\
		\cite{li2021short} & Cluster & Withd. & 1 Hour & - & Spatial-Temporal Memory Network \\
		\cite{li2023improving} & Cluster & Withd. & 1 Hour & - & Irregular Convolutional LSTM \\
		\cite{li2015traffic} & Cluster & Withd. and return & 1 Hour & - & GB \\
		\cite{pan2019predicting} & Cluster & Withd. and return & 1 Hour & 24 Hours & LSTM, Deep Neural Networks \\
		\cite{zhang2018short} & System & Demand & 15 Min. & - & LSTM \\
		\cite{cantelmo2020low} & System & Demand & 1 Day & - & Low Dimensionality Model \\
		  \cite{li2019citywide} & Station, Cluster, System & Withd. and return & 1 Hour & - & Hierarchical Consistency Prediction Model \\
	\end{tabular}
	\end{adjustbox}
	\caption{Subset of the most relevant contributions in the field of forecasting for BSS}
	\label{tab:ref}
\end{table}

\myparagraph{Forecasting and relocation}

Different approaches have been discussed in the literature when the forecast of user behavior and decision-making for BSSs are both addressed. 
These vary both in the approach used to 
forecast user behavior (e.g., supervised algorithm), as well as in the way the decision-making that happens on the basis of this information is tackled (e.g., optimization, simulation). 
A summary of the contributions is presented in Table \ref{tab:litrev_relevant} reporting, for each of the selected references, the forecasting strategy and, for the decision-making side, the approach. The table indicates also whether decisions are defined in a static or dynamic setting, and whether the decision involves routing constraints for the relocating vehicles. In the following, we will describe the main contributions addressing both aspects, highlighting those that investigate how the former is functional to the latter. 

One of the earliest contributions addressing both tasks is \cite{caggiani2013dynamic}. The authors present a dynamic simulator where, given the withdrawal demand for each station, the destination choice is simulated to assess the arrival time of the users. Demand for withdrawal and return is estimated with two Neural Networks (NN). Bike relocation is addressed by means of a non-linear optimization problem to find the pairwise flow of bikes among stations minimizing relocation costs and number of lost users. Routing constraints are considered assuming one vehicle is available for the relocation of bikes. 
Results are provided for a BSS with 5 stations in three different days and with three different levels of demand for each day, but the impact of the forecast on the quality of the relocation solution is not evaluated.
\cite{o2015data} consider both the overnight and mid-rush relocation, discussing both as static optimization problems. The demand for bikes in each station is computed as a historical average. As historical data of bike trips differ from trip demand due to censoring, for each time window past data are only considered when a station is not empty. Accounting for high traffic, mid-rush relocation is performed by small vehicles carrying a limited number of bikes between pairs of closely-located stations. Overnight, relocation is tackled by planning truckloads of bikes in order to maximize the number of stations in which relocation takes place. Results are reported on 50 real-life instances extracted from the June 2013 data of the Citibike system, and compare the performance of a mixed-integer formulation with those of a greedy algorithm.
\cite{alvarez2016optimizing} use statistical inference to estimate the parameters of the non-homogeneous Poisson process representing the withdrawals and returns considering 1 hour time intervals for each station in each type of day (i.e., working day, Saturdays, and Sundays).
Unsatisfied demand is forecasted via an approximated Monte Carlo simulation. The forecasts are used to solve a static bike relocation problem with routing constraints. Results are reported for the BSS of Palma de Mallorca, including 28 stations, specifically for the historical data spanning one week of November 2013. 
\cite{liu2016rebalancing} tackle the forecasting task by means of a $k$-Nearest Neighbor-based regression model for the forecasting of the withdrawals in each station. An inter-station bike transit model is used to forecast the station returns and return time after each withdrawal event. A non-linear formulation is proposed to solve a static bike relocation problem. A clustering method is used to group stations which are then assigned to one vehicle in the presented solution algorithm. While the forecasting algorithm is shown to have better performance than the compared models, the effect of different forecasts on the relocation problem is not discussed.
\cite{mrazovic2018deep} propose the use of a multi-output Long Short-Term Memory (LSTM) to tackle the task of forecasting the withdrawal and return requests for the following 8 half-hours. The multiple steps ahead capability of the presented model is used to identify the target inventory maximizing the number of consecutive time intervals during which the station remains balanced.

\cite{fan2019distributed} use Random Forest (RF) to forecast the target level of each station for the next day. The relocation problem is modeled as a static vehicle routing problem with pickup and delivery, with an ant colony optimization heuristic being presented for its solution. The performance of the presented heuristic is compared with that of a branch \& cut implementation for the layout of different BSSs, while the performance of the methodological approach as a whole is discussed for the BSS active in Hangzhou, China.
A similar static relocation problem is considered in \cite{yoshida2019practical} where the issue of forecasting bike usage and relocating bikes to mitigate imbalances is tackled in three steps. First, the number of withdrawals and returns for the next day is forecasted by means of a graph convolutional network. Second, an optimization problem is solved to define the optimal inventory level for each station. Third, based on the solution of the previous step, a static vehicle routing problem is solved to plan the relocation of bikes. 
\cite{lee2020optimal} present a bike relocation strategy that considers both daily demand forecasting, achieved by means of a boosting algorithm, and bike relocation with vehicle routing constraints.
In the routing problem defined for bike relocation, travel times are time-dependent—varying between peak and off-peak hours—and the fleet size used for relocation is assumed to be flexible, allowing the number of relocating vehicles to adapt so that all operations can be completed within a prescribed time limit.
A genetic algorithm is proposed for the solution of the problem and tested on the real data of November 2018 of the 95 stations of the Gangnam district of Seoul, South Korea.
\cite{chiariotti2020bike} present an optimization framework to evaluate the effect of the joint use of bike relocation and user incentives on the overall service quality of a BSS. The process of departures and arrivals at each station is modeled as an independent Markov-modulated Poisson process whose arrival and departure rates follow historic patterns. This description of the system is used to define the expected amount of time a station is either empty or full. Both static and dynamic relocation strategies are considered. While the capacity of the relocating vehicles is considered, the rebalancing operations and the traveling of the vehicles between stations is assumed to be instantaneous. 
In \cite{cho2021efficiency}, multiple RFs are used to forecast user behavior. Most importantly, different models are used to forecast the use of the system by casual and regular users. A linear formulation is proposed for the static bike relocation problem, also considering routing constraints. Different repositioning strategies are tested, 
based either on safety stock buffers or on forecast–error thresholds used to determine whether a station requires service.
Computational results are used to identify the conditions under which each strategy performs best, depending on demand patterns and forecast accuracy.
\cite{cipriano2021data} present a pattern mining approach between neighboring stations to detect stations with critical occupancy levels. The rules embed the information of whether a station is positively or negatively critical (i.e., lack of empty stands or lack of bikes, respectively) and are used to define stations where bikes are loaded from or unloaded to when performing bike relocation. 
The effect of a number of hyperparameters of the pattern mining technique on the relocation operations is reported. 
Given the unsupervised nature of the forecasting technique, no result is reported on how the quality of the forecast helps the decision-making process. Instead, results are reported discussing the effect of the number of relocating vehicles and of the data partitioning strategy used for rule extraction, when relocation happens twice a day (i.e., at 6:00 and 15:00).
\cite{lin2022demand} proposes to exploit demand patterns to cluster stations with similar demand trends. A non-linear demand scaling technique is presented to represent how stations with higher capacity can afford greater demands for withdrawals and returns.
The scaled demand is then used to build a forecast with a rule based on exponential decay. Lower and upper values are dynamically defined for the forecast value to trigger the relocation of bikes. 
In \cite{huang2022monte} a density-based clustering technique is used to group stations based on the usage pattern and relocating vehicles are assigned to each cluster. The demand for withdrawals and returns is forecasted by averaging historical data.
The relocation problem for each vehicle is solved dynamically by means of a Monte Carlo tree search. The results of the relocation approach are compared to those obtained when no relocation is performed and to those of alternative approaches based on greedy algorithms and reinforcement learning. 

Interesting considerations on the effect of forecasting on the downstream decision-making task are drawn in \cite{regue2014proactive} and \cite{gammelli2022predictive}.
\cite{regue2014proactive} present a four step approach. First, a station-level demand forecasting model is defined, based on Gradient Boosting (GB). Second, an inventory level model is defined where each station is represented as a queuing model to determine whether the initial inventory of a given station is enough to serve the demand of a certain future time period.
Third, a stochastic linear integer program is used to define the number of bikes to be picked up or delivered to each station. Finally, a single-vehicle routing problem is solved to plan the operations each time it returns to the depot. Each vehicle trip is set to last a maximum of 20 minutes, which also coincides with the time unit of the forecasting task.
The impact of GB on relocation performance is compared to those of Linear Regression (LR) and of the random uniform sampling of a target within station capacity. 
\cite{gammelli2022predictive} present a recurrent NN-based approach to forecast withdrawal and return rates for one station in a station-based BSS, which are used to define the parameters of a non-homogeneous Poisson process. An inventory model using the forecasted rates is defined to determine the target level of a station at the beginning of a day. The target for each station is  to minimize a user dissatisfaction function, defined by the integral of penalties for lost withdrawals and returns over time.
It is worth noting that no consideration is presented on the operational logistics required to achieve such target values, such as the planning or execution of bike relocation within the system.
The performance of the presented forecasting model compared to other more classical approaches is shown both with respect to the forecasting error and in the effect on inventory optimization. Regarding the latter, results obtained with the presented forecasting method are shown to be 40\% closer to oracle performance compared to models such as LR. Interestingly, the authors highlight that the proposed forecasting methodology performs best in the forecasting task when withdrawals and returns are forecasted jointly but the best results for inventory optimization are obtained when the two quantities are forecasted separately. The authors describe how this result is due to the fact that the optimal decision is influenced by the cumulative difference between pickups and returns, rather than only their separate evolution over the day.

\begin{table}[ht]
\centering
\begin{adjustbox}{max size={1\textwidth}{1\textheight}}
\begin{tabular}{l|llll}
\multicolumn{2}{l}{} & \multicolumn{3}{c}{Decision-making} \\ \cline{3-5} 
Reference & Forecasting & Approach	& Dynamic	& Routing \\ \hline
\cite{caggiani2013dynamic} & Supervised & Simulation Optimization & Yes & Yes \\ 
\cite{mrazovic2018deep} & Supervised & Rule-based & No & No \\
\cite{gammelli2022predictive} & Supervised & Optimization & No & No \\
\cite{liu2016rebalancing} & Supervised & Optimization & No & Yes \\
\cite{fan2019distributed} & Supervised & Optimization & No & Yes \\
\cite{yoshida2019practical} & Supervised & Optimization & No & Yes \\
\cite{lee2020optimal} & Supervised & Optimization & No & Yes \\
\cite{cho2021efficiency} & Supervised & Optimization & No & Yes \\
\cite{regue2014proactive} & Supervised & Optimization & Yes & Yes \\
\cite{tomaras2018modeling} & Smoother & Optimization & Yes & Yes \\ 
\cite{o2015data} & Historical average & Optimization & No & Yes \\
\cite{huang2022monte} & Historical average & Monte Carlo-based & No & No \\
\cite{alvarez2016optimizing} & Statistical inference & Optimization & No & Yes \\
\cite{cipriano2021data} & Pattern recognition & Rule-based & Yes & No \\
\cite{chiariotti2020bike} & Markov Process & Optimization & Yes & Partly \\
\cite{lin2022demand} & Historical data & Rule-based & Yes & No \\
\end{tabular}
\end{adjustbox}
\caption{Literature addressing both the task of forecasting and that of decision-making.}
\label{tab:litrev_relevant}
\end{table}

\section{Forecasting approaches and models}\label{sec:forecast}

In this section, we first describe the process of creating the dataset, starting from the raw data collected by a BSS-operating company. Then, we introduce the ML algorithms used and, finally, we discuss the two approaches employed to leverage the features available in the dataset.

Companies operating BSSs maintain records of the system usage, typically in a database that logs where, when, and by whom each bike is picked up and dropped off. The raw data from these logs have been processed to obtain, for each station, the count of bike withdrawals and returns in predefined time intervals.
To balance between accuracy and data sparsity while still allowing to plan for bike relocation, user operations are grouped into half-hour intervals, counting withdrawals and returns per interval and calculating the net demand.
Using a longer time unit for the interval (e.g., one hour) reduces the number of intervals in a dataset spanning several years, which can facilitate the identification of temporal dependencies over longer periods. However, it can obscure intra-interval variability.
For example, if, within a single time unit, both a withdrawal and a return are logged, the net demand might appear the same regardless of the actual sequence of the events. 
Yet, if a return precedes a withdrawal, the actual number of bikes in the station is different than that of the reverse scenario. 
Defining a short time unit (e.g., 1 minute, 5 minutes) mitigates this issue but increases the number of intervals over the same time period (e.g., one year), potentially complicating the identification of meaningful temporal patterns in the data. 
Moreover, a shorter time interval often increases data sparsity, making demand forecasting more challenging.

In this study a half-hour time unit is used while considering a forecasting horizon of one year. 
This enables the forecasts to be used in a station-based BSS simulator to guide strategic and tactical decision-making while still representing the dynamic relocation of bikes during the simulated day. Furthermore, it allows a detailed evaluation of the quality of the forecast of each ML approach as well as an assessment of the impact that each approach has on the simulation of the BSS and on the decision-making for the system. The explanatory features resulting from transforming the raw data provided by the company are of different nature and are listed in below. 
For each feature, the type is also reported, with ``numerical'' indicating features that take continuous values, ``categorical'' indicating features that take discrete values identifying categories, and ``boolean'' indicating features that have a $\{0, 1\}$ domain.

\begin{itemize}
    \setlength\itemsep{0em}
	\item Half-hour within the day (categorical)
	\item Day of month (categorical)
	\item Day of week (categorical)
	\item Month of the year (categorical)
	\item Progressive week index (numerical)
	\item Station identifier (categorical)
\end{itemize}

These features have been integrated with the following meteorological and date-related information.

\begin{itemize}
    \setlength\itemsep{0em}
	\item Average daily temperature (numerical)
	\item Average wind speed (numerical)
	\item Public holiday (boolean)
	\item School day (boolean)
	\item Occurrence of fog during the day (boolean)
	\item Occurrence of rain during the day (boolean)
	\item Occurrence of snow during the day (boolean)
	\item Occurrence of storm during the day (boolean)
\end{itemize}

Days have been classified as ``school days'' if they are not main holidays. For the case of Bicimia presented in Section \ref{sec:bsmob}, these are summer holidays (typically from the second week of June to the first week of September) and winter holidays (typically, the week before and the one after new year's eve). The definition of weather-related features as the characteristics of a day, rather than of an half-hour time slot, allows the users of the simulator to more easily define their values over a long horizon while also allowing the model to learn from the available meteorological information.

Different ML-based regression methods have been tested to forecast the net demand. As shown in Section \ref{sec:litrev}, the review of the literature highlighted how the most capable techniques for this task are frequently based on Neural Networks or defined as ensemble models. This informed the methodological choices of the present study, leading to the adoption of Neural Networks, Long Short-Term Memory networks, Random Forest, and Light Gradient-Boosting Machine (LightGBM).

Neural Networks (NN) are a classical ML method. They consist of interconnected layers of neurons that process input data to recognize patterns and make forecasts. Each neuron applies an activation function to a weighted sum of its inputs and passes the result. This output is used as an input by the neurons of subsequent layers or as the output of the model. 
The NN architecture typically includes an input layer, one or more hidden layers, and an output layer. In this study, a feedforward architecture is used, where information moves in one direction from input to output.
Their ability to learn complex, non-linear relationships made Neural Networks a staple of machine learning.

Long Short-Term Memory (LSTM) networks \citep{hochreiter1997long} are a specialized class of Recurrent Neural Networks (RNN) specifically designed to model sequential and time-dependent data. Unlike feedforward Neural Networks, LSTM networks process sequences of input observations and explicitly account for temporal dependencies by maintaining an internal state that evolves over time. This makes them particularly suitable for forecasting tasks where past observations influence future values, such as time series prediction. 
The core component of an LSTM network is the memory cell, which is regulated by a set of gating mechanisms: the input gate, forget gate, and output gate. These gates control the flow of information into, within, and out of the cell, allowing the network to retain relevant information while discarding irrelevant or outdated one, mitigating the vanishing gradient problem that affects standard RNNs. Regularization techniques such as L2 regularization and dropout are frequently used during the training of both NNs and LSTM networks to mitigate the risk of overfitting.

Random Forest \citep{breiman2001random} is a classical ensemble learning method that combines multiple decision trees through a technique called bagging (bootstrap aggregation). In bagging, each tree is built on a random subset of the data, created by sampling with replacement, and a random subset of features. For regression purposes, the final prediction is made by averaging the predictions of all trees.
This technique reduces overfitting and improves generalization, making Random Forest particularly robust to the presence of noise and outliers in the data.

LightGBM~\citep{ke2017lightgbm} is a framework developed with the aim of improving upon existing gradient boosting methods.
Gradient boosting (see \cite{friedman2001greedy}) is an ensemble technique that builds a model through a sequence of weak learners (i.e., supervised learning models with limited complexity), where each new weak learner is trained to correct the residuals (errors) of the previous one, as weighted by a learning rate.
LightGBM builds trees with a depth-first rather than breadth-first approach, selecting the node in the tree with the highest loss. This can lead to better accuracy and faster training time when compared to traditional boosting frameworks.
LightGBM also exploits two novel techniques called Gradient-based One-Side Sampling and Exclusive Feature Bundling. The first focuses on the most informative data points by keeping instances with large gradients, with the aim of improving training speed by selecting a subset of data points that contribute most to the information gain. Exclusive Feature Bundling reduces the dimensionality by bundling mutually exclusive categorical features. 
This is especially useful when dealing with high-dimensional sparse data, as it allows a more compact representation of features without losing important information.

Two approaches have been tested to exploit the available data with the considered ML algorithms.
In the first, the station identifier is used as an explanatory feature, with one algorithm being trained on data regarding all stations. We will refer to this as the ``global'' variant.
The second is a station-specific approach that trains one model for each station of the system. We will refer to this variant as ``local''. This alternative approach was formulated in response to the distinct patterns of bike station usage observed among users in the real-world case at hand, which will be described in depth in Section \ref{sec:bsmob} (see Figure \ref{fig:net_stations}).
Note that given that LSTM models are inherently designed to learn from sequential data, a temporal structure is imposed on the bike-sharing dataset. The LSTM-global model is therefore designed to perform multivariate forecasting of the vector of net demand, with one component of the vector for each station.

To help the assessment of the quality of each model, results have been reported for an additional model, termed \histShift{} (\HS{}), defined as the historical data shifted by 52 weeks backward (or, more precisely, by 17472 half-hours). In this model, the prediction of the net demand of bikes
for a station in a half-hour interval is defined as the observed net demand for that station in that same interval 52 weeks prior. This model was included to assess the importance of capturing the temporal characteristics of net demand of bike, besides its magnitude. 

\section{Simulation framework} \label{sec:simfw}

\begin{figure}[H]
	\centering
	\includegraphics[width=0.7\linewidth]{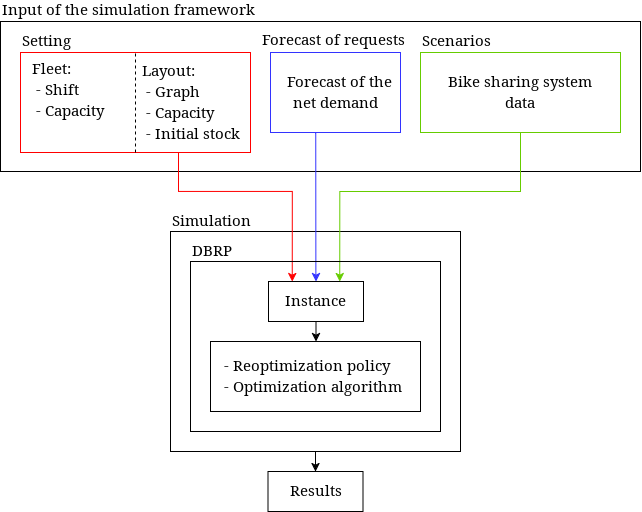}
	\caption{The structure of the  modified simulation framework.}
	\label{fig:simulation}
\end{figure}

To analyze the benefit of long term forecasting for strategic decision-making BSSs, a simulation framework has been adopted in this study,  based on the one presented in \cite{angelelli2022simulation}. 
In this section, we describe the main components of the simulator and illustrate the changes we introduced for the present contribution. The modified scheme of the simulator is depicted in Figure \ref{fig:simulation}.

\begin{figure}[H]
	\centering
    \includegraphics[width=0.7\linewidth]{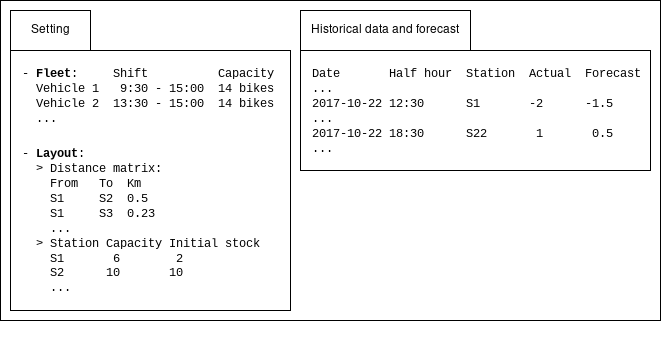}
	\caption{An example of the data defining an instance.}
	\label{fig:instance}
\end{figure}

The simulator requires several inputs, described in the following. The first set of inputs outlines the operational environment of the system, referred to as \textit{setting}, and is defined by two components, referred to as the \textit{layout} and the \textit{fleet}. The layout is the data describing the stations in terms of their location, capacity, and initial stock at the beginning of the simulated day. It also includes the information regarding the traveling time and distance between each pair of stations, denoted in the figure as ``Graph''. The fleet is the data describing the vehicles for the bike relocation operations. This information includes the number of vehicles and, for each vehicle, its capacity and the start and end time of its shift in the day being simulated.

For each simulated day, additional inputs are required regarding the requests of the users, namely, the way requests manifest throughout the day in the BSS, which is referred to as a \textit{scenario}, and the \textit{forecast of the requests} used to plan the relocation.
In \cite{angelelli2022simulation}, scenarios are generated through the definition of non-homogeneous Poisson processes for the returns and withdrawals of bikes of each station. The definition of these processes is dependent on a number of characteristics of the simulated day which also shape the forecasting model currently defined by the company operating the Bicimia BSS in Brescia, Italy (see Section \ref{sec:cmexpl}).
In this paper, actual historical data will be used to provide scenarios for the simulation allowing us to draw insights that are more applicable to the challenges faced by company operating the BSS. Scenarios will therefore be defined by actual actions of users in each tested day, as recorded by the BSS. An example of the data defining an instance in the simulator is reported in Figure \ref{fig:instance}.
The simulation framework will be used to evaluate the benefit of long term ML-based forecasting for strategic decisions in BSSs. To this end, the different forecasting  methods presented in Section \ref{sec:forecast} will be tested to generate the forecast of the net demand for the case study of the Bicimia BSS, and will be compared to the strategy applied by the company operating the BSS. The tested methods will be fed with values for the features defined in Section \ref{sec:forecast}, hence providing predictions for the net demand for each station in the required half-hour.

The definition of the setting, of the scenario and of the forecast of the requests allows the simulation framework to  represent the behavior of the users of the system and the bike relocation happening during the day. The latter is defined through the solution of a Dynamic Bike Relocation Problem (see \cite{shui2020review}), with a reoptimization strategy triggered by key events (e.g., vehicle arrival or departure from stations).

\section{A case study: Bicimia} \label{sec:bsmob}

\subsection{The data}
\begin{figure}[H]
	\centering	\includegraphics[width=0.4\linewidth]{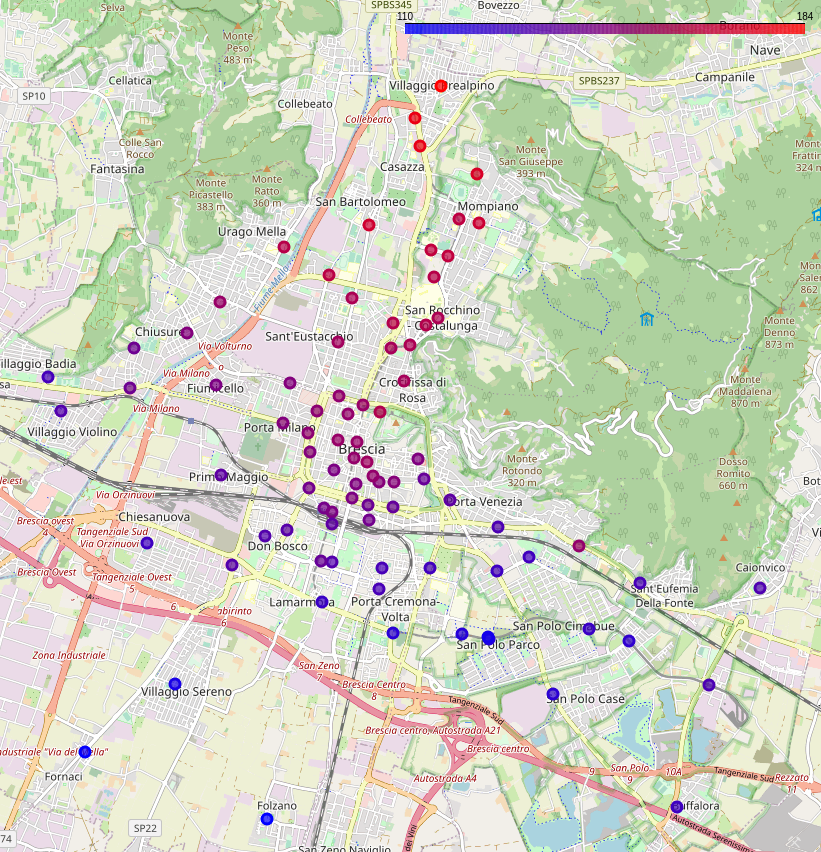}
	\caption{The geographical layout and elevation of the stations of the Bicimia service.}
	\label{fig:map}
\end{figure}

This contribution builds upon a collaboration of the authors with Brescia Mobilità (see \cite{angelelli2022simulation}), running a BSS called Bicimia which serves the municipality of Brescia, Italy, and a few neighboring towns. To date, the system counts 94 stations, covering an area of approximately 55 $km^2$ (13600
acres). While the city is located at a meeting point between the Alps and the Po valley, the urban portion of the municipality of Brescia is rather flat, with a 64 meter elevation difference between the highest and the lowest station. The geographical layout and elevation of the stations are highlighted in Figure \ref{fig:map}. 
At the time of writing, the capacity of the stations ranges from 5 to 26 bikes, with a capacity of 10 bikes being the most common. 

\begin{figure}[H]
	\centering
	\includegraphics[width=0.4\linewidth]{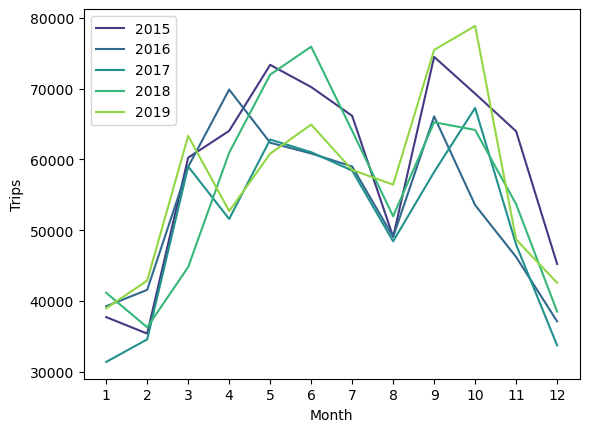}
	\caption{The monthly number of trips by year.}
	\label{fig:years_by_month}
\end{figure}

\begin{figure}[H]
	\centering
	\includegraphics[width=0.4\linewidth]{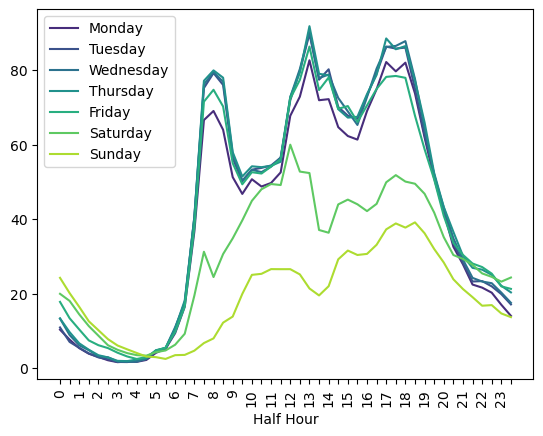}
	\caption{The average number of trips for each half-hour interval by day of the week.}
	\label{fig:days_by_hh}
\end{figure}

\begin{figure}[H]
	\centering
	\begin{subfigure}{0.32\textwidth}
		\centering
		\includegraphics[width=\linewidth]{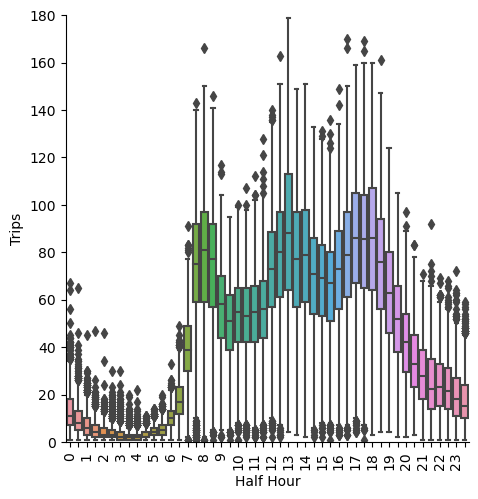}
		\caption{Working days}
	\end{subfigure}
	\begin{subfigure}{0.32\textwidth}
		\centering
		\includegraphics[width=\linewidth]{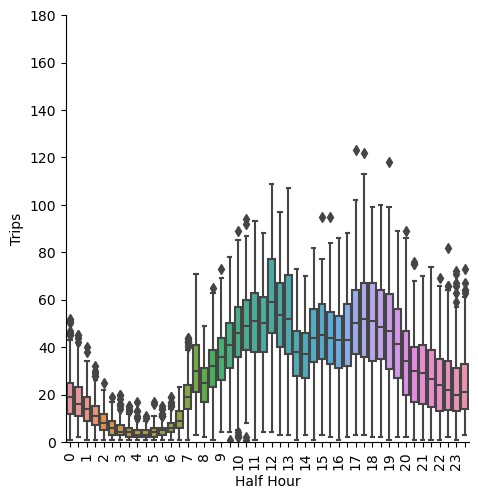}
		\caption{Saturday}
	\end{subfigure}
	\begin{subfigure}{0.32\textwidth}
		\centering
		\includegraphics[width=\linewidth]{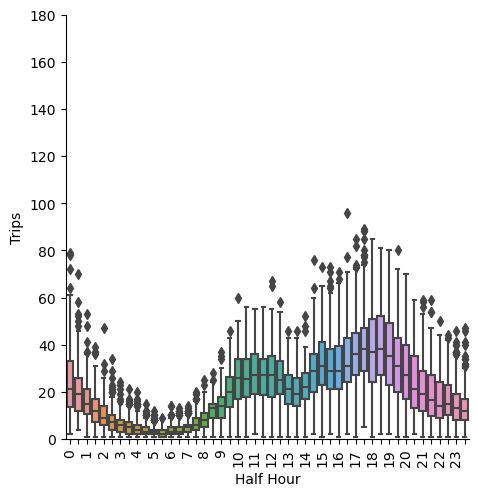}
		\caption{Sunday}
	\end{subfigure}
	\caption{Half-hour usage distribution for working days, Saturday, and Sunday.}
	\label{fig:bo_days_by_hh}
\end{figure}

Data from the company were made available to the authors for the years in the interval $[2015,2019]$ in the form of the logs of bike withdrawals and returns. As can be observed from Figure \ref{fig:years_by_month}, the usage of the system has been rather stable during this time interval, with no major trends across the years. 
The average usage for 30 minutes intervals is reported in Figure \ref{fig:days_by_hh}. Akin to other bike sharing systems (e.g., \cite{alvarez2016optimizing}), 
working days exhibit a similar pattern, whereas Saturday and Sunday show different trends. Slight decreases can be observed on Mondays and in Friday's morning and evening peaks. 
The company has attributed the former to the absence of relocation operations on Sundays,
causing a less than ideal distribution of bikes to satisfy the users' transportation needs. The latter could be caused by a subset of commuters selecting Friday as the day of the week to work from home. The distribution of system usage by half-hour, as defined by the trip starting time, is reported in Figure \ref{fig:bo_days_by_hh}.
Differently from most BBSs (e.g., \cite{kim2019graph,boufidis2020development,cantelmo2020low}), the Bicimia BSS exhibits a ``noon'' peak during working days.
This peak is also observed during Saturdays and is explained by the intense use of the system by high school students.
Finally, it can be observed that, as the week progresses, the use of the system during night time appears to shift to later hours.

\begin{figure}[H]
	\centering
	\begin{subfigure}{0.49\textwidth}
		\centering
		\includegraphics[width=0.7\linewidth]{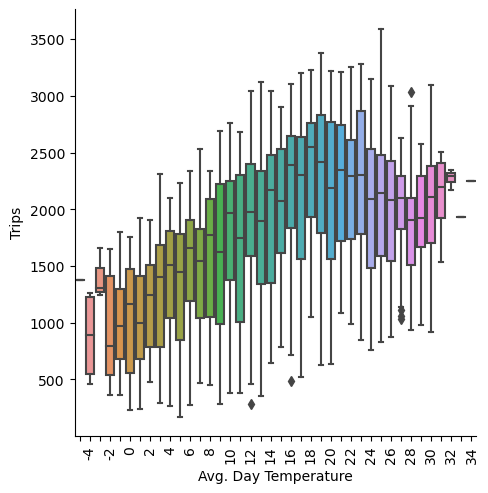}
		\caption{Daily number of trips w.r.t. average day temperature.}
	\end{subfigure}
	\begin{subfigure}{0.49\textwidth}
		\centering
		\includegraphics[width=0.7\linewidth]{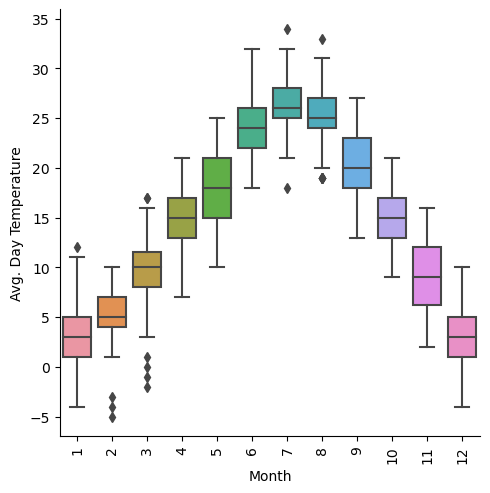}
		\caption{Average day temperature by month in Brescia in the considered period.}
	\end{subfigure}
	\caption{The effect of temperature on the usage of the BSS and the average temperatures in the considered period.}
	\label{fig:temperature}
\end{figure}

\begin{figure}[H]
	\centering
	\begin{subfigure}{0.24\textwidth}
		\centering
		\includegraphics[width=\linewidth]{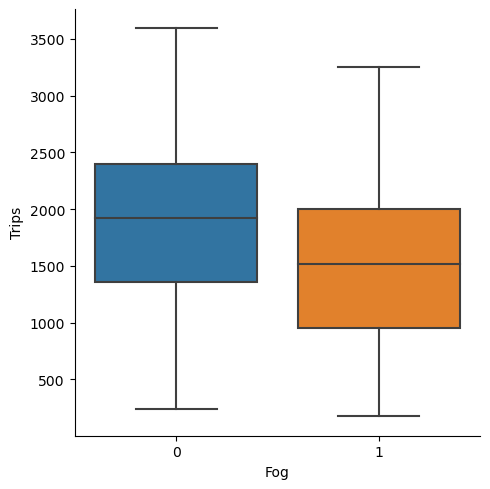}
		\caption{Fog}
	\end{subfigure}
	\begin{subfigure}{0.24\textwidth}
		\centering
		\includegraphics[width=\linewidth]{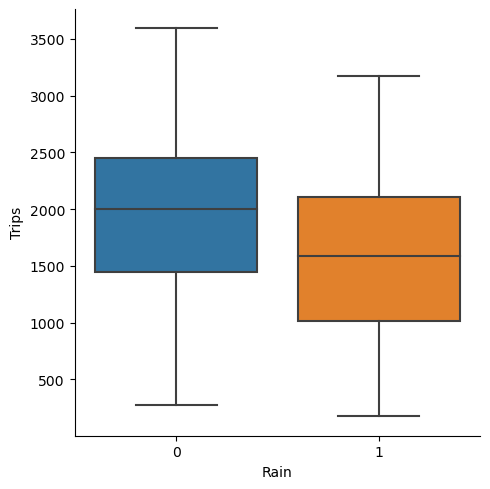}
		\caption{Rain}
	\end{subfigure}
	\begin{subfigure}{0.24\textwidth}
		\centering
		\includegraphics[width=\linewidth]{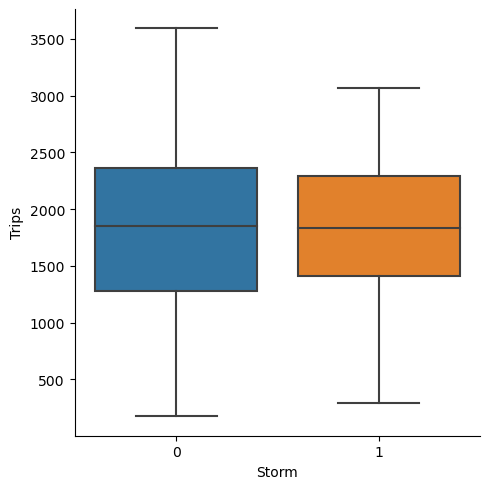}
		\caption{Storm}
	\end{subfigure}
	\begin{subfigure}{0.24\textwidth}
		\centering
		\includegraphics[width=\linewidth]{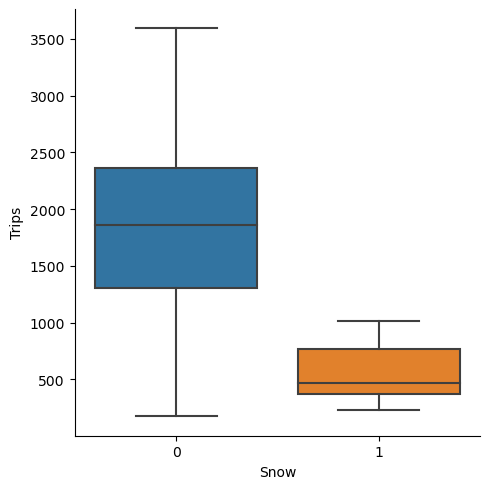}
		\caption{Snow}
	\end{subfigure}
	\caption{The effect of weather phenomena on the use of the system.}
	\label{fig:precip}
\end{figure}

\begin{figure}[H]
	\centering
	\includegraphics[width=0.4\linewidth]{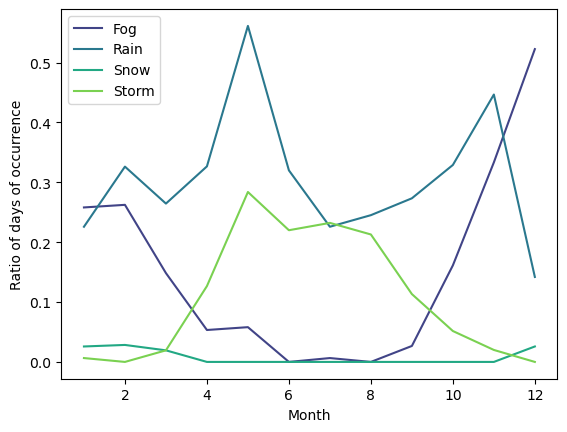}
	\caption{The monthly occurrence of each phenomenon as ratio of the days for each month.}
	\label{fig:phenomena_by_month}
\end{figure}

Data from the company have been integrated with daily historical meteorological data. This allowed an investigation into the effect of different weather characteristics on system utilization. 
In particular, the impact of daily average temperature on bike usage and the monthly distribution of the temperature in Brescia are reported in Figure \ref{fig:temperature}, while the impact of different weather phenomena is shown in Figure \ref{fig:precip}. 
As expected, the daily number of trips has a peak with mildly warm temperature in the interval [15°, 25°], and decreases with more extreme temperatures, especially colder ones. As reported in Section \ref{sec:forecast}, with respect to weather phenomena, the information available is a boolean variable on the occurrence of each phenomenon during that day, that is, whether the phenomenon has occurred during that day or not.
As can be observed, fog and rain appear to have a similar effect on the usage of the system. Storms appear to have little effect on users habits, probably due to their short duration. Finally, snow is shown to drastically reduce bike use. This is partially due to its occurrence in colder temperatures. The monthly occurrence of each phenomenon is reported in Figure \ref{fig:phenomena_by_month} as ratio of the days of occurrence, for each month.

\begin{figure}[H]
	\centering
	\begin{subfigure}{0.48\textwidth}
		\centering
		\includegraphics[width=0.7\linewidth]{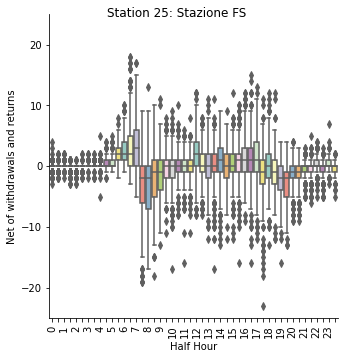}
		\caption{BSS station outside the Railway station.}
	\end{subfigure}
	\begin{subfigure}{0.48\textwidth}
		\centering
		\includegraphics[width=0.7\linewidth]{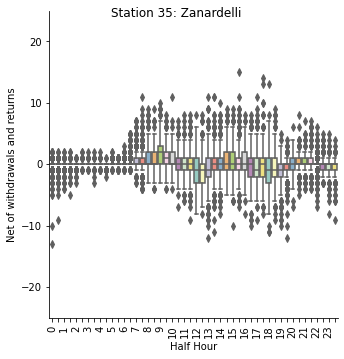}
		\caption{City center, popular spot BSS station.}
	\end{subfigure}
	\begin{subfigure}{0.48\textwidth}
		\centering
		\includegraphics[width=0.7\linewidth]{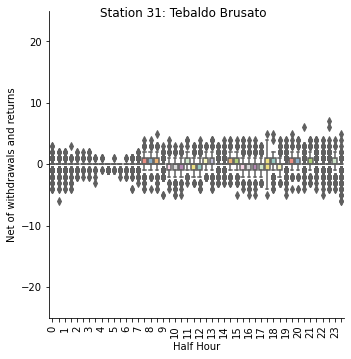}
		\caption{City center, lightly used BSS station.}
	\end{subfigure}
	\begin{subfigure}{0.48\textwidth}
		\centering
		\includegraphics[width=0.7\linewidth]{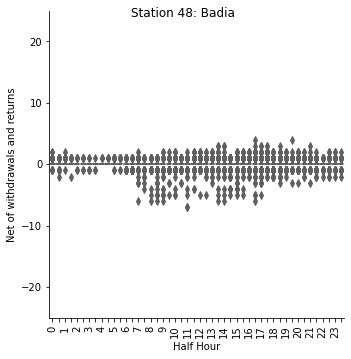}
		\caption{Suburb BSS station.}
	\end{subfigure}
	\caption{The trend of the net demand of bikes in different stations.}
	\label{fig:net_stations}
\end{figure}

A crucial aspect of bike relocation is the net between bike returns and withdrawals, taking positive value when the returns are greater than withdrawals and negative value when the opposite is true. It can be easily seen how this quantity manifests the imbalance in user demand in a station. To showcase the different behavior of the net across the system, in Figure \ref{fig:net_stations} values have been reported in half-hour intervals for different stations in the system. Different behaviors can be observed for the stations. The station with highest fluctuations is the one near the city's railway station (Figure \ref{fig:net_stations}(a)), where fluctuations caused by outbound and inbound commuters can be observed in the morning, in the [5:30, 7:30] and [8:00, 9:30] intervals, respectively. In the evening, a downward peak is shown in the time slot [18:00, 20:00], possibly caused by outbound commuters returning home after work. Figure \ref{fig:net_stations} also shows the behavior for other stations in the city center, one being a popular spot in the city and the other being a lightly used station (Figure \ref{fig:net_stations}(b) and (c), respectively). A station in the suburbs is also shown (Figure \ref{fig:net_stations}(d)).

Currently, the company operates a fleet of three vehicles, each with a capacity of 14 bikes. Operators are scheduled as follows. During working days (i.e., Monday to Friday), two operators work the [7:00 - 15:00] shift, while one works the [11:30 -19:30] shift. During Saturdays, one operator works the [7:00 - 13:00] shift. Relocation is not performed on Sundays.

\subsection{The current forecasting model} \label{sec:cmexpl}

Currently, the company does not use ML techniques to forecast withdrawals and returns of bikes at the stations. For this reason, in \cite{angelelli2022simulation} the characteristics of a day influencing the forecast of the requests have been identified according to the experience of the operators. Such characteristics are the following:
\begin{itemize}
	\item type of day, i.e., working day, Saturday, or Sunday;
	\item month of the year;
	\item weather, i.e., sunny or rainy.
\end{itemize}

For each combination of these characteristics, a day in 2019 has been identified to provide the profile of withdrawals and returns to be taken as a reference for each station. The resulting forecasting model, therefore, assumes that the users will behave accordingly to the profile of the selected day taken as reference. Henceforth, we will refer to this model as \companyModel{} (\CM).

\section{Results}\label{sec:results}

\subsection{Forecasting}

\begin{figure}[H]
	\centering	\includegraphics[width=0.5\linewidth]{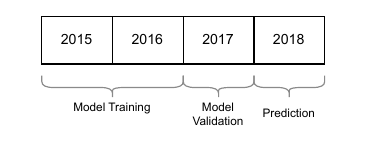}
	\caption{Usage of the data for the different years.}
	\label{fig:years}
\end{figure}

As defined in Section \ref{sec:cmexpl}, the \companyModel{} is defined based on historical data drawn from the year 2019. Hence, to provide a fair comparison of the performance of the models, experiments have been carried out on the observations of 2015 and 2016, model validation and light hyper-parameter tuning on the observations of 2017, leaving 2018 as unobserved data for all the approaches (see Figure \ref{fig:years}). 
It is worth noting that this data handling approach also facilitates the retraining of the forecasting models as new data become available. This is particularly important since the evolution of the system could cause data drift and need to periodically re-evaluate strategic and tactical decisions.

For all ML methods, a light hyper-parameter tuning has been performed on the validation set. Specifically, for Neural Networks tuning was applied on the number of hidden layers, on the number of nodes in each of the hidden layers, and on the strength of the L2 regularization term. 
For LSTM, tuning was applied to  the number of LSTM layers, the number of nodes in each of the layers, as well as the number of nodes in the feedforward (dense) layer following the LSTM ones, and the dropout percentage.
For Random Forest, tuning was performed on the number of estimators (decision trees) and on their maximum depth. For LightGBM, tuning was performed on the number of estimators (decision trees) and their maximum depth, as well as on the learning rate, and on the fraction of samples and features used to train each tree.

The performance of each model is reported in Table \ref{tab:MAEMSE}, for both 2017 and 2018, in terms of Mean Squared Error (MSE) as defined in the equation: 
$$MSE = \frac{1}{n} \sum_{i=1}^n (y_i - \hat{y}_i)^2,$$
where $y_i$ denotes the true net demand of bikes and $\hat{y}_i$ the forecast of the model, for each observation $i$ in the dataset. For each year, the ``\% gap'' row reports the percentage gap with respect to the best performing model.

\begin{table}[ht]
\centering
\begin{adjustbox}{max size={0.9\textwidth}{0.9\textheight}}
\begin{tabular}{cccccccccccc}
& \multicolumn{1}{c}{} & & \multicolumn{4}{c}{Global} & \multicolumn{4}{c}{Local}	\\ \cmidrule(lr){4-7} \cmidrule(lr){8-11} 
& \HS & \CM & NN & LSTM & RF & LightGBM & NN & LSTM & RF & LightGBM  \\ \cmidrule(lr){2-2} \cmidrule(lr){3-3} \cmidrule(lr){4-7} \cmidrule(lr){8-11}
2017 & 2.090 & 2.124  & 1.118 & 1.031 & 1.095 & 1.112  & 1.053 & 1.041 & 1.056 & \textbf{1.031} \\
\% gap & 102.65 & 105.98 & 8.43 & 0.03 & 6.23 & 7.86 & 2.11 & 0.96 & 2.36 & - \\ \cmidrule(lr){2-11}
2018 & 2.023 & 2.115 & 1.124 & 1.079 & 1.111 & 1.120  & 1.092 & 1.077 & 1.089 & \textbf{1.067} \\
\% gap & 89.59 & 98.21 & 5.28 & 1.07 & 4.13 & 4.91   & 2.30 & 0.93 & 2.03 & - \\ \cmidrule(lr){2-11}
\end{tabular}
\end{adjustbox}
\caption{The MSE of the tested models for 2017 and 2018.}
\label{tab:MAEMSE}
\end{table}

As can be observed, for both years, the best performance is obtained with the LightGBM-local model with the closest competitor being LSTM-global for the validation data and LSTM-local for the test data.
This highlighting the capability of gradient boosting to extract meaningful relationships between explanatory features and target variable. This capability is further discussed in Section \ref{sec:featimpo}, where the feature importance output of LightGBM-local is examined.
Notably, on the 2018 data, for all tested methods the local approach is shown to perform better than the global counterpart.
This result can be justified by the different behavior of the users in different stations (see Figure \ref{fig:net_stations}), which renders beneficial the adoption of models that are tailored to capture the behavior of users in each station. The result is consistent with those found in the literature (e.g., \cite{regue2014proactive}). 

Interestingly, \companyModel{} (\CM{}) and \histShift{} (\HS{})
are shown to be the worst performing models in both years, with both having about twice the MSE than LightGBM-local.
Finally, we note that while the performance of all ML-based models is shown to slightly deteriorate from 2017 to 2018, the opposite is true for \HS{} and \CM{}.
This could be due to the fact that, to compute predictions for 2018, data from 2017 and 2019 are used by \HS{} and \CM{}, respectively, that is, data with at most one year of lag with the prediction, whereas ML-based models are limited to the use of data up to 2016, that is, data with at least two years of lag.
Moreover, the behavior of the ML-based models is not surprising, as these have been validated on the observations of 2017 while 2018 remained unobserved.

\begin{figure}[H]
	\centering
	\includegraphics[width=0.4\linewidth]{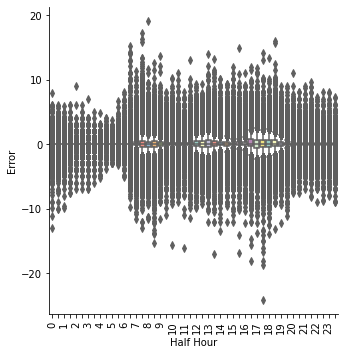}
	\caption{Distribution of the error by half-hour.}
	\label{fig:error_halfhour}
\end{figure}

\begin{figure}[H]
	\centering
	\begin{subfigure}{0.32\textwidth}
		\centering
		\includegraphics[width=\linewidth]{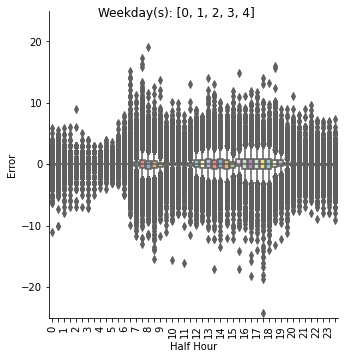}
		\caption{}
	\end{subfigure}
	\begin{subfigure}{0.32\textwidth}
		\centering
		\includegraphics[width=\linewidth]{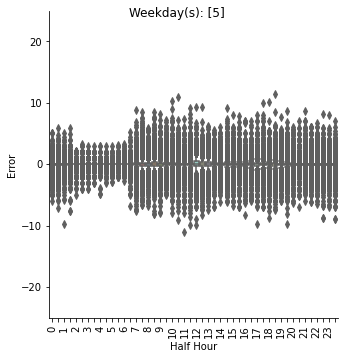}
		\caption{}
	\end{subfigure}
	\begin{subfigure}{0.32\textwidth}
		\centering
		\includegraphics[width=\linewidth]{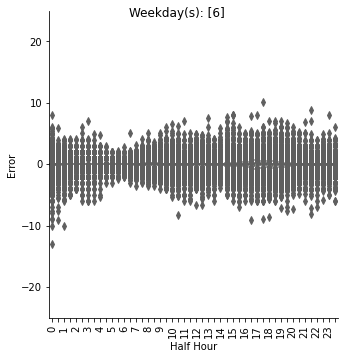}
		\caption{}
	\end{subfigure}
	\caption{Distribution of the error by day of the week: (a) working days, (b) Saturdays, (c) Sundays.}
	\label{fig:error_weekday}
\end{figure}

\begin{figure}[H]
	\centering
	\begin{subfigure}{0.48\textwidth}
		\centering
		\includegraphics[width=0.7\linewidth]{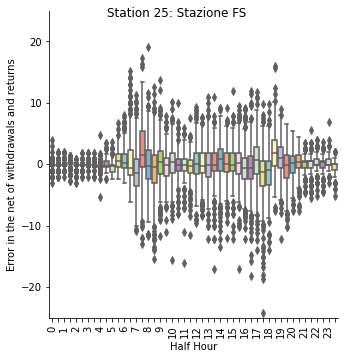}
		\caption{Railway station}
	\end{subfigure}
	\begin{subfigure}{0.48\textwidth}
		\centering
		\includegraphics[width=0.7\linewidth]{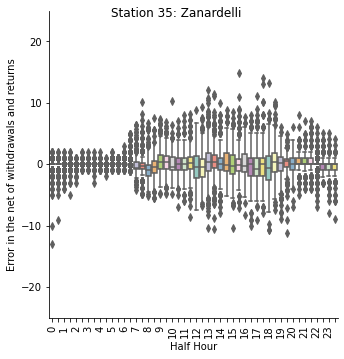}
		\caption{City center, popular spot station.}
	\end{subfigure}
	\begin{subfigure}{0.48\textwidth}
		\centering
		\includegraphics[width=0.7\linewidth]{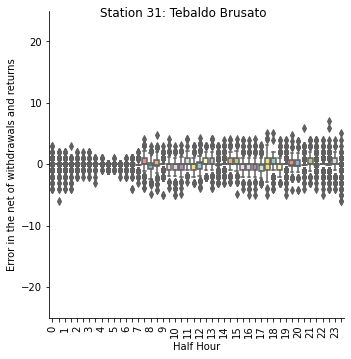}
		\caption{City center, lightly used station.}
	\end{subfigure}
	\begin{subfigure}{0.48\textwidth}
		\centering
		\includegraphics[width=0.7\linewidth]{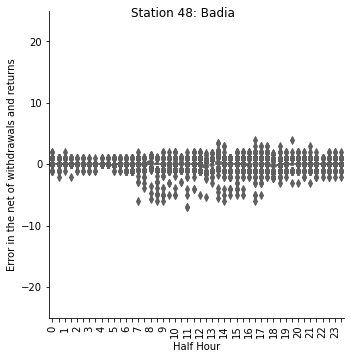}
		\caption{Suburb station.}
	\end{subfigure}
	\caption{The trend of the error in the stations reported in Figure \ref{fig:net_stations}.}
	\label{fig:error_stations}
\end{figure}

The distribution of the error of the LightGBM-local approach for all stations is reported in Figures \ref{fig:error_halfhour}, \ref{fig:error_weekday}, and \ref{fig:error_stations}, according to various explanatory features, and is calculated for each observation of the respective half-hour across the 365 observations of 2018.
As can be seen from Figure \ref{fig:error_halfhour}, which shows the distribution of the error for each half-hour of the day, the error of the model appears to be rather constant 
during daytime, with the exception of a reduction in variance with the decrease in the usage of the system during nighttime. A slight increase can be observed in the error variance and asymmetry during the morning and evening peak, further indicating an increase in the difficulty of the forecasting task as the system is used the most.
This interpretation is validated by the results shown in Figure \ref{fig:error_weekday}, where the same plot is depicted for working days, Saturdays, and Sundays, as well as Figure \ref{fig:error_stations}, showing how higher errors are observed in stations with higher traffic, particularly, the one located close to the railway station and the one located in a popular spot of the city center.

\subsection{Feature importance}\label{sec:featimpo}

\begin{figure}[H]
	\centering
	\includegraphics[width=0.5\linewidth]{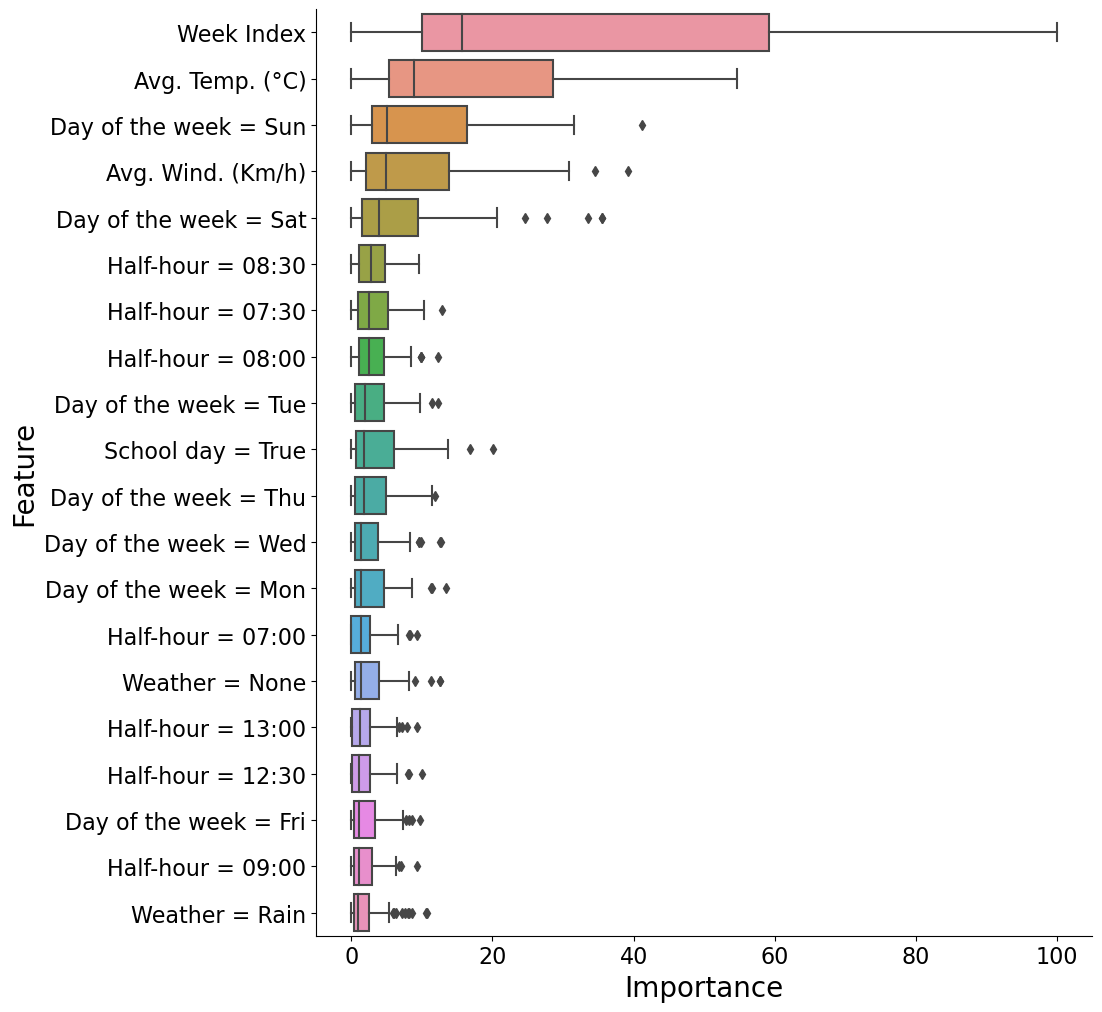}
	\caption{Boxplot of the feature importance for the 20 most important features across all models.}
	\label{fig:featimpo}
\end{figure}

To study the relationship between the explanatory features and the target variable for the LightGBM-local model, we resorted to the native capability of LightGBM to provide feature importance as an output. In particular, by default, LightGBM defines the importance of a feature by the number of times it is used in the model, i.e., the number of times a split is made in a decision tree based on the feature.
Given the one model per station approach of LightGBM-local, feature importance has been analyzed as the normalized distribution of the result of the model of each station, reported in Figure \ref{fig:featimpo}.

First and foremost, from Figure \ref{fig:featimpo}, we note the importance of the week index, which shows how temporal trends play a role of primary importance in forecasting the demand of bikes in the system. The relevant role of weather on the forecast of the net use of the system is also observed, with temperature and wind being, on average, the second and fourth most important features. 
As expected, temporal variables  have a clear role. 
Day of the week is also shown to provide important information on bike usage, 
with the related variables appearing
among the 20 most important features. This is especially evident for weekend days, which exhibit different dynamics in user behavior, as shown in Figures \ref{fig:days_by_hh} and \ref{fig:bo_days_by_hh}.
With respect to half-hour information, half-hours pertaining to the morning peak, i.e., from 7:00 to 9:00, are shown to have the highest importance. Interestingly, and denoting the peculiarity of the Bicimia BSS, half-hours pertaining to the noon peak are shown to have higher importance than those of the evening peak. This result, jointly with the relevance of the features reporting whether one day is a school day and whether the observation is in September (the first month of school), highlights the considerable use of the system by students.

\subsection{Simulation}

The forecasts of each of the models discussed in the previous section have been used in the simulation framework to assess the quality of the resulting bike relocation, as measured by the total number of missed withdrawals, returns, and their sum, for each day of year 2018. 
The computations of missed withdrawals and returns are carried out according to the formulas reported in \cite{angelelli2024one}. In short, each time a station is empty and a withdrawal request occurs, the number of missed withdrawals increases. Similarly, each time a station is full and a return request occurs, the number of missed returns increases.
Results are reported in Table \ref{tab:simres} as averages over the days when operators are scheduled to perform bike relocation, i.e., Monday to Saturday, for a total of 313 days.
The table also reports the average daily number of kilometers (Total Km) driven by the vehicles and the average daily number of relocated bikes, with the former being defined as the sum of the distances between the stations visited by each vehicle, and the latter being the average of a counter variable that keeps track of the loading and unloading operations between each vehicle and a station.

\begin{table}[H]
\centering
\begin{adjustbox}{max size={0.9\textwidth}{0.9\textheight}}
\begin{tabular}{rrrrrrrrrrrrr}
	& \multicolumn{2}{c}{} & \multicolumn{4}{c}{Global} & \multicolumn{4}{c}{Local}	\\ \cmidrule(lr){4-7} \cmidrule(lr){8-11}
	& \HS & \CM & NN & LSTM & RF & LightGBM & NN & LSTM & RF & LightGBM \\ \cmidrule(lr){2-2} \cmidrule(lr){3-3} \cmidrule(lr){4-7} \cmidrule(lr){8-11}
Missed withdrawals & 61.29 & 60.65 & \textbf{50.00} & 56.34 & 56.70 & 56.77 & 60.22 & 55.94 & 54.48 & 53.62 \\ 
Missed returns & \textbf{44.13} & 44.60 & 57.14 & 54.18 & 49.80 & 51.53 & 48.80 & 44.45 & 46.84 & 45.92 \\ \cmidrule(lr){2-11}
Total missed requests & 105.42 & 105.25 & 107.13 & 110.51 & 106.50 & 108.30 & 109.02 & 100.39 & 101.33 & \textbf{99.54} \\ \cmidrule(lr){2-11}
\begin{tabular}[c]{@{}c@{}} \% gap \\ from \CM \end{tabular} & 0.16 & - & 1.79 & 5.00 & 1.18 & 2.90 & 3.58 & -4.61 & -3.72 & -5.43 \\ \cmidrule(lr){2-11}
Total Km & 252.22 & 254.55 & 236.65 & 220.39 & 222.27 & 219.47 & 245.23 & 247.77 & 234.78 & 246.26 \\
Relocated bikes & 336.17 & 332.69 & 343.85 & 312.21 & 307.52 & 313.32 & 307.21 & 317.58 & 318.78 & 322.52 \\
\end{tabular}
\end{adjustbox}
\caption{The average daily number of missed withdrawals, returns and requests, the percentage gap from the \companyModel{}, the total distance (in km), and the number of relocated bikes.}
\label{tab:simres}
\end{table}

Among the tested models, and consistently with the forecasting performance, the best results are achieved with the LightGBM-local approach, providing an average improvement of about 5.43\% with respect to the \companyModel. Two other local models, namely LSTM and RF, yield measurable benefits in this regard.
The improvement achieved by the use of the forecasts of the LightGBM-local is even more consistent when considering that the average amount of requests that are missed when perfect information is used in the simulation is 62.21 in 2018.
If we subtract this non-reducible amount of missed requests to both, the improvement achieved by the LightGBM-local approach becomes 13.28\%.
This result is achieved while also reducing the total distance traveled by the vehicles, with an average reduction of more than 8 Km per day, and less work for the operators, reducing the average daily number of relocated bikes by about ten units.

Overall, the difficulty of the forecasting task, characterized by a long forecasting horizon with a small temporal unit definition, and of the subsequent bike relocation is reflected in the results obtained by the simulation of the system, with the \companyModel{}
achieving comparable results to those of some of the ML models.
This also highlights how different approaches in the definition and the selection of the model impact the performance of the bike sharing system as defined in the simulation framework.
In particular, the function used to measure the system performance also plays an important role in establishing which method guarantees the best results. For instance, in the Bicimia BSS users are given locks which provide an alternative way to return bikes when no stands are available at a station (besides looking for an alternative station with an empty stand). Conversely, no alternative exists for the lack of bikes at a station, e.g., bike reservation, making the minimization of missed withdrawals a possible objective in a practical implementation of this approach. As shown in Table \ref{tab:simres}, with such objective all ML-based approaches outperform the \CM{} and NN-global proves to be the best model.

\begin{figure*}[ht]
	\centering
		\centering 
		\includegraphics[width=0.5\textwidth]{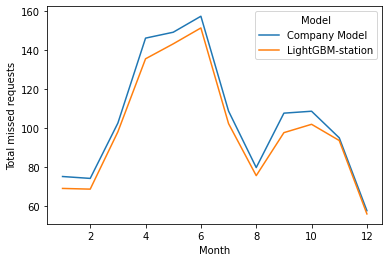}
    \caption{Distribution of the average daily missed requests by month.} 
	\label{fig:sim_by_month}
\end{figure*}

To analyze the improvement provided by LightGBM-local over the \companyModel{}, Figure \ref{fig:sim_by_month} shows the performance of the simulator, by month, as measured by the total missed requests for 2018.
LightGBM-local is shown to provide the greatest advantage over the \companyModel{} in the months of April to June, September, and October, and the least advantage in the months of November and December. This result could be explained by the use of the system by workers and, especially, students going to school, which is rather regular in school months with mild temperature. These users are possibly the most regular users of the system, with a predefined and somewhat static daily schedule. The decrease of the use of the system by these users might be the cause of an increased difficulty of the forecasting of the net level of withdrawals and returns.

\begin{figure*}[ht]
		\centering \includegraphics[width=0.5\textwidth]{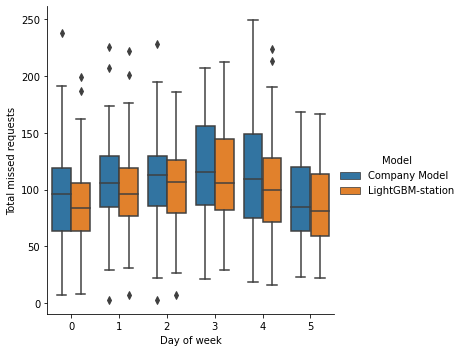}
    \caption{Distribution of the simulator performance by day of the week in terms of the distribution of total missed requests.}
	\label{fig:sim_by_dayofweek} 
\end{figure*}

Figure \ref{fig:sim_by_dayofweek}
shows the performance of the simulator 
by day of the week, in terms of the distribution of total missed requests, when provided with the \CM{} and LightGBM-local forecasts. Note that Sundays are not reported as no operator is scheduled to work in this day of the week. 
As can be observed, both approaches show a worsening of the performance for the relocation as the week goes by, especially until Thursday. 
Our intuition is that the current fleet of relocating vehicles is not appropriately sized for the Bicimia BSS. Possibly, relocating vehicles are not able to keep up with the imbalances created during the working days, which are characterized by a higher number of requests. This situation worsens as the week progresses and is only mitigated on Saturdays, in which a lower demand is observed.

To showcase the impact of different forecasting approaches on strategic and tactical decision-making, in Tables \ref{tab:chiari} and \ref{tab:LIGHTGBM}, the performance of the simulator, when either \CM{} or LightGBM-local are used, are reported for different combinations of vehicles used in the morning and afternoon shifts. In the tables, the lowest average number of missed requests for each value of the fleet size is highlighted in bold.
These tables provide a valid guidance to make and reassess decisions about the sizing of the fleet and the scheduling of the vehicles. In this regard, the bold values in the tables highlight the best combination of vehicles in the morning and afternoon shift for a fleet of a given size. For instance, in Table \ref{tab:chiari}, 99.48 is the lowest average number of missed requests for a fleet of four vehicles, achieved when three vehicles are scheduled in the morning and one in the afternoon. As can be observed by comparing the bold values of Tables \ref{tab:chiari} and \ref{tab:LIGHTGBM}, 
the best scheduling of the vehicles for each given fleet size might differ between \CM{} and LightGBM-local, i.e., for a fleet of three and five vehicles, highlighting how the two approaches could lead to different strategic and tactical decisions. Furthermore, the tables show how, for a  number of vehicles from 1 to 5, the best result is always achieved when LightGBM-local is used.

\begin{table}[H]
	\centering
	\begin{tabular}{c|r|rrrrrr}
		\multicolumn{1}{r}{} & \multicolumn{1}{r}{} & \multicolumn{6}{c}{\# Vehicles - Afternoon shift} \\
		\cmidrule{3-8}	\multicolumn{1}{r}{} & \multicolumn{1}{r}{} & \multicolumn{1}{c}{0}	 & \multicolumn{1}{c}{1}	 & \multicolumn{1}{c}{2}	 & \multicolumn{1}{c}{3}	 & \multicolumn{1}{c}{4}	 & \multicolumn{1}{c}{5} \\
		\cmidrule{3-8}	\multirow{6}[1]{*}{\begin{tabular}[l]{@{}c@{}} \# Vehicles - \\ Morning shift \end{tabular}} & 
		0 & \textbf{145.69} & 144.77 & 138.36 & 132.68 & 129.11 & 127.52  \\ 
		& 1 & \textbf{131.69} & 124.91 & 116.85 & 112.67 & 110.72  & \\ 
		& 2 & \textbf{116.16} & 108.43 & 102.84 & 99.80 &  & \\ 
		& 3 & \textbf{106.12} & \textbf{99.48} & \textbf{94.58} &  & & \\ 
		& 4 & 101.52 & 95.22 &  & & & \\ 
		& 5 & 99.46 &  & & & & \\ 
			\end{tabular}
	\caption{Average daily number of missed requests for each combination of morning and afternoon shifts, when \CM{} is used.}
	\label{tab:chiari}
 \end{table}

 \begin{table}[H]
	\centering
	\begin{tabular}{c|r|rrrrrr}
		\multicolumn{1}{r}{} & \multicolumn{1}{r}{} & \multicolumn{6}{c}{\# Vehicles - Afternoon shift} \\
		\cmidrule{3-8}	\multicolumn{1}{r}{} & \multicolumn{1}{r}{} & \multicolumn{1}{c}{0}	 & \multicolumn{1}{c}{1}	 & \multicolumn{1}{c}{2}	 & \multicolumn{1}{c}{3}	 & \multicolumn{1}{c}{4}	 & \multicolumn{1}{c}{5}	 \\
		\cmidrule{3-8}	\multirow{6}[1]{*}{\begin{tabular}[l]{@{}c@{}} \# Vehicles - \\ Morning shift \end{tabular}} & 
		0 & \textbf{145.69} & 136.43 & 131.04 & 128.69 & 127.39 & 126.75  \\ 
		& 1 & \textbf{122.15} & 112.56 & 108.93 & 107.29 & 106.30  & \\ 
		& 2 & \textbf{110.39} & \textbf{102.29} & 99.02 & 97.34 &  & \\ 
		& 3 & 103.74 & \textbf{96.43} & 93.57  & & & \\ 
		& 4 & 100.04 & \textbf{92.52} & & & & \\ 
		& 5 & 97.57 
			\end{tabular}
	\caption{Average daily number of missed requests for each combination of morning and afternoon shifts, when LightGBM-local is used.}
	\label{tab:LIGHTGBM}
 \end{table}

\section{Conclusions} \label{sec:conclusions}

In this paper, the impact of ML-based forecasting on tactical and strategical decision-making for a station-based BSS has been analyzed.
Different machine learning methods have been tested in order to forecast a year long interval of the net demand of bikes
for each station of a bike sharing system, with an half-hour time unit. This particular setting of the forecasting task is required in order to generate forecasts which are defined on a time unit small enough for operational decisions concerning the relocation to be embedded in the simulation, and on a horizon  which is long enough to assess the impact of strategic and tactical decisions.The results show how machine learning forecasting may provide a moderate advantage over the rule-based approach currently in use by the company, with a decrease in the number of missed requests of bike withdrawals and returns of around 5.4\%.

A major takeaway of this study is the difficulty of  long-term forecasting for strategic and tactical decision-making, especially when the simulation of daily operational decisions is required to assess their impact.
Moreover, the performance of the forecasting models and of the simulation  was shown to be impacted by the approach, either global or local, taken in the forecasting models. This leads to the conclusion that much attention should be put in how to use the BSS station information in the forecasting model.

Finally, the selection of the best forecasting model was shown to be dependent on the metric used to define the quality of the downstream simulation, highlighting the importance of testing the forecasting model with respect to its impact on the selected performance measure of the system.

The present study can be extended in several directions. 

Further investigating the relationship between forecasting and BSS simulation is an area of future research. 
In this work, user behavior and demand  have been represented, for both the forecast and the simulation, by the net demand of bikes
at each station.
More detailed representations, such as, demand for a trip from a specific bike station to another, or from the actual origin to destination of the user, would allow for a more accurate estimation of the relationship between request forecasting and decision-making for bike sharing systems. 

\section*{Acknowledgments}

The authors are grateful to three anonymous reviewers whose reports helped them substantially improve a previous version of this paper.
They would also like to express their gratitude to Brescia Mobilità SpA, with special thanks to Massimo Chiari, for the collaboration that formed the foundation of this research.

This research of M. Grazia Speranza has been supported and funded by the European Union (EU) and Italian Ministry for Universities and Research (MUR), National Recovery and Resilience Plan (NRRP), within the project “Sustainable Mobility Center (MOST)”, 2022–2026, CUP D83C22000690001, Spoke N° 7, “CCAM, Connected networks and Smart Infrastructures” and the grant PRIN 2022-NRRP.

\section*{Disclosure statement}
The authors report there are no competing interests to declare.

\bibliographystyle{abbrvnat}
\bibliography{biblio}

@article{mor2026optimization,
   title={Optimization and supervised learning for decision making: competitors or partners?},
    author={Mor, Andrea and Orsenigo, Carlotta and Speranza, M. Grazia},
    journal={to appear on Journal of the Operational Research Society},
    year={2026},
    doi={10.1080/01605682.2026.2616419}
}

@article{hochreiter1997long,
  title={Long short-term memory},
  author={Hochreiter, Sepp and Schmidhuber, J{\"u}rgen},
  journal={Neural computation},
  volume={9},
  pages={1735--1780},
  year={1997},
  publisher={MIT press}
}

@misc{meddin2024,
  author = {O’Brien, Oliver and DeMaio, Paul and Rabello, Renata and Chou, Steve and Benicchio, Thiago},
  title = {{The Meddin Bike-sharing World Map Report
2024}},
  howpublished = "\url{https://bikesharingworldmap.com/reports/bswm\_2024report.pdf}",
  year = {2024}, 
  note = "[Online; accessed 28-11-2025]"
}

@article{li2019citywide,
  title={Citywide bike usage prediction in a bike-sharing system},
  author={Li, Yexin and Zheng, Yu},
  journal={IEEE Transactions on Knowledge and Data Engineering},
  volume={32},
  number={6},
  pages={1079--1091},
  year={2019},
  publisher={IEEE}
}

@article{ma2022short,
  title={Short-term prediction of bike-sharing demand using multi-source data: a spatial-temporal graph attentional {LSTM} approach},
  author={Ma, Xinwei and Yin, Yurui and Jin, Yuchuan and He, Mingjia and Zhu, Minqing},
  journal={Applied Sciences},
  volume={12},
  number={3},
  pages={1161},
  year={2022},
  publisher={MDPI}
}

@article{luo2021predicting,
  title={Predicting travel demand of a docked bikesharing system based on {LSGC-LSTM} networks},
  author={Luo, Jing and Zhou, Dai and Han, Zhaolong and Xiao, Guangnian and Tan, Yunlong},
  journal={IEEE Access},
  volume={9},
  pages={92189--92203},
  year={2021},
  publisher={IEEE}
}

@article{torres2024forecasting,
  title={Forecasting the usage of bike-sharing systems through machine learning techniques to foster sustainable urban mobility},
  author={Torres, Jaume and Jim{\'e}nez-Mero{\~n}o, Enrique and Soriguera, Francesc},
  journal={Sustainability},
  volume={16},
  number={16},
  pages={6910},
  year={2024},
  publisher={MDPI}
}

@article{teusch2023systematic,
  title={A systematic literature review on machine learning in shared mobility},
  author={Teusch, Julian and Gremmel, Jan Niklas and Koetsier, Christian and Johora, Fatema Tuj and Sester, Monika and Woisetschl{\"a}ger, David M and M{\"u}ller, J{\"o}rg P},
  journal={IEEE Open Journal of Intelligent Transportation Systems},
  volume={4},
  pages={870--899},
  year={2023},
  publisher={IEEE}
}

@article{zhang2023relocation,
  title={Relocation-related problems in vehicle sharing systems: A literature review},
  author={Zhang, Ruiyou and Kan, Haiyu and Wang, Zhaoming and Liu, Zhujun},
  journal={Computers \& Industrial Engineering},
  volume={183},
  pages={109504},
  year={2023},
  publisher={Elsevier}
}

@article{angelelli2024one,
title = {The one-station bike repositioning problem},
journal = {Discrete Applied Mathematics},
volume = {357},
pages = {173-196},
year = {2024},
author = {Angelelli, Enrico and Mor, Andrea and Speranza, M. Grazia},
}

@article{haider2018inventory,
  title={Inventory rebalancing through pricing in public bike sharing systems},
  author={Haider, Zulqarnain and Nikolaev, Alexander and Kang, Jee Eun and Kwon, Changhyun},
  journal={European Journal of Operational Research},
  volume={270},
  pages={103--117},
  year={2018},
  publisher={Elsevier}
}

@article{ghosh2017dynamic,
  title={Dynamic repositioning to reduce lost demand in bike sharing systems},
  author={Ghosh, Supriyo and Varakantham, Pradeep and Adulyasak, Yossiri and Jaillet, Patrick},
  journal={Journal of Artificial Intelligence Research},
  volume={58},
  pages={387--430},
  year={2017}
}

@article{raviv2013static,
  title={Static repositioning in a bike-sharing system: models and solution approaches},
  author={Raviv, Tal and Tzur, Michal and Forma, Iris A},
  journal={EURO Journal on Transportation and Logistics},
  volume={2},
  pages={187--229},
  year={2013},
  publisher={Elsevier}
}

@article{friedman2001greedy,
  title={Greedy function approximation: a gradient boosting machine},
  author={Friedman, Jerome H},
  journal={Annals of Statistics},
  pages={1189--1232},
  year={2001},
  publisher={JSTOR}
}

@article{breiman2001random,
  title={Random forests},
  author={Breiman, Leo},
  journal={Machine Learning},
  volume={45},
  pages={5--32},
  year={2001},
  publisher={Springer}
}

@article{ke2017lightgbm,
  title={Lightgbm: A highly efficient gradient boosting decision tree},
  author={Ke, Guolin and Meng, Qi and Finley, Thomas and Wang, Taifeng and Chen, Wei and Ma, Weidong and Ye, Qiwei and Liu, Tie-Yan},
  journal={Advances in Neural Information Processing Systems},
  volume={30},
  year={2017}
}

@inproceedings{yoshida2019practical,
  title={Practical end-to-end repositioning algorithm for managing bike-sharing system},
  author={Yoshida, Akihiro and Yatsushiro, Yosuke and Hata, Nozomi and Higurashi, Tatsuru and Tateiwa, Nariaki and Wakamatsu, Takashi and Tanaka, Akira and Nagamatsu, Kenichi and Fujisawa, Katsuki},
  booktitle={2019 IEEE International Conference on Big Data (Big Data)},
  pages={1251--1258},
  year={2019},
  organization={IEEE}
}

@article{huang2022monte,
  title={Monte {C}arlo tree search for dynamic bike repositioning in bike-sharing systems},
  author={Huang, Jianbin and Tan, Qinglin and Li, He and Li, Ao and Huang, Longji},
  journal={Applied Intelligence},
  volume={52},
  pages={1--16},
  year={2022},
  publisher={Springer}
}

@inproceedings{tomaras2018modeling,
  title={Modeling and predicting bike demand in large city situations},
  author={Tomaras, Dimitrios and Boutsis, Ioannis and Kalogeraki, Vana},
  booktitle={2018 IEEE International Conference on Pervasive Computing and Communications (PerCom)},
  pages={1--10},
  year={2018},
  organization={IEEE}
}

@article{fan2019distributed,
  title={Distributed forecasting and ant colony optimization for the bike-sharing rebalancing problem with unserved demands},
  author={Fan, Yiwei and Wang, Gang and Lu, Xiaoling and Wang, Gaobin},
  journal={PloS One},
  volume={14},
  pages={e0226204},
  year={2019},
  publisher={Public Library of Science San Francisco, CA USA}
}

@article{lin2022demand,
  title={A demand-centric repositioning strategy for bike-sharing systems},
  author={Lin, Ying-Chih},
  journal={Sensors},
  volume={22},
  pages={5580},
  year={2022},
  publisher={MDPI}
}

@article{cipriano2021data,
  title={A data-driven based dynamic rebalancing methodology for bike sharing systems},
  author={Cipriano, Marco and Colomba, Luca and Garza, Paolo},
  journal={Applied Sciences},
  volume={11},
  pages={6967},
  year={2021},
  publisher={MDPI}
}

@article{cho2021efficiency,
  title={Efficiency comparison of public bike-sharing repositioning strategies based on predicted demand patterns},
  author={Cho, Jung-Hoon and Seo, Young-Hyun and Kim, Dong-Kyu},
  journal={Transportation Research Record},
  volume={2675},
  pages={104--118},
  year={2021},
  publisher={SAGE Publications Sage CA: Los Angeles, CA}
}

@inproceedings{liu2016rebalancing,
  title={Rebalancing bike sharing systems: A multi-source data smart optimization},
  author={Liu, Junming and Sun, Leilei and Chen, Weiwei and Xiong, Hui},
  booktitle={Proceedings of the 22nd ACM SIGKDD International Conference on Knowledge Discovery and Data Mining},
  pages={1005--1014},
  year={2016}
}

@article{regue2014proactive,
  title={Proactive vehicle routing with inferred demand to solve the bikesharing rebalancing problem},
  author={Regue, Robert and Recker, Will},
  journal={Transportation Research Part E: Logistics and Transportation Review},
  volume={72},
  pages={192--209},
  year={2014},
  publisher={Elsevier}
}

@article{alvarez2016optimizing,
  title={Optimizing the level of service quality of a bike-sharing system},
  author={Alvarez-Valdes, Ramon and Belenguer, Jose M and Benavent, Enrique and Bermudez, Jose D and Mu{\~n}oz, Facundo and Vercher, Enriqueta and Verdejo, Francisco},
  journal={Omega},
  volume={62},
  pages={163--175},
  year={2016},
  publisher={Elsevier}
}

@article{lee2020optimal,
  title={Optimal relocation strategy for public bike system with selective pick-up and delivery},
  author={Lee, Euntak and Son, Bongsoo and Han, Youngjun},
  journal={Transportation Research Record},
  volume={2674},
  pages={325--336},
  year={2020},
  publisher={SAGE Publications Sage CA: Los Angeles, CA}
}

@inproceedings{o2015data,
  title={Data analysis and optimization for (citi) bike sharing},
  author={O'Mahony, Eoin and Shmoys, David},
  booktitle={Proceedings of the AAAI Conference on Artificial Intelligence},
  volume={29},
  year={2015}
}

@inproceedings{li2015traffic,
  title={Traffic prediction in a bike-sharing system},
  author={Li, Yexin and Zheng, Yu and Zhang, Huichu and Chen, Lei},
  booktitle={Proceedings of the 23rd SIGSPATIAL International Conference on Advances in Geographic Information Systems},
  pages={1--10},
  year={2015}
}

@article{caggiani2021toward,
  title={Toward sustainability: Bike-sharing systems design, simulation and management},
  author={Caggiani, Leonardo and Camporeale, Rosalia},
  journal={Sustainability},
  volume={13},
  pages={7519},
  year={2021},
  publisher={MDPI}
}

@inproceedings{zhang2018short,
  title={Short-term prediction of bike-sharing usage considering public transport: A {LSTM} approach},
  author={Zhang, Cheng and Zhang, Linan and Liu, Yangdong and Yang, Xiaoguang},
  booktitle={2018 21st International Conference on Intelligent Transportation Systems (ITSC)},
  pages={1564--1571},
  year={2018},
  organization={IEEE}
}

@article{li2023improving,
  title={Improving short-term bike sharing demand forecast through an irregular convolutional neural network},
  author={Li, Xinyu and Xu, Yang and Zhang, Xiaohu and Shi, Wenzhong and Yue, Yang and Li, Qingquan},
  journal={Transportation Research Part C: Emerging Technologies},
  volume={147},
  pages={103984},
  year={2023},
  publisher={Elsevier}
}

@article{li2021short,
  title={Short-term forecast of bicycle usage in bike sharing systems: a spatial-temporal memory network},
  author={Li, Xinyu and Xu, Yang and Chen, Qi and Wang, Lei and Zhang, Xiaohu and Shi, Wenzhong},
  journal={IEEE Transactions on Intelligent Transportation Systems},
  volume={23},
  pages={10923--10934},
  year={2021},
  publisher={IEEE}
}

@article{collini2021deep,
  title={Deep learning for short-term prediction of available bikes on bike-sharing stations},
  author={Collini, Enrico and Nesi, Paolo and Pantaleo, Gianni},
  journal={IEEE Access},
  volume={9},
  pages={124337--124347},
  year={2021},
  publisher={IEEE}
}

@article{wang2018short,
  title={Short-term prediction for bike-sharing service using machine learning},
  author={Wang, Bo and Kim, Inhi},
  journal={Transportation Research Procedia},
  volume={34},
  pages={171--178},
  year={2018},
  publisher={Elsevier}
}

@article{ermagun2018urban,
  title={Urban trails and demand response to weather variations},
  author={Ermagun, Alireza and Lindsey, Greg and Loh, Tracy Hadden},
  journal={Transportation Research Part D: Transport and Environment},
  volume={63},
  pages={404--420},
  year={2018},
  publisher={Elsevier}
}

@article{sohrabi2020real,
  title={Real-time prediction of public bike sharing system demand using generalized extreme value count model},
  author={Sohrabi, Soheil and Paleti, Rajesh and Balan, Lacramioara and Cetin, Mecit},
  journal={Transportation Research Part A: Policy and Practice},
  volume={133},
  pages={325--336},
  year={2020},
  publisher={Elsevier}
}

@article{zhou2019reliable,
  title={A reliable traffic prediction approach for bike-sharing system by exploiting rich information with temporal link prediction strategy},
  author={Zhou, Yan and Li, Yanxi and Zhu, Qing and Chen, Fen and Shao, Junming and Luo, Yunxin and Zhang, Yeting and Zhang, Pengcheng and Yang, Weijun},
  journal={Transactions in GIS},
  volume={23},
  pages={1125--1151},
  year={2019},
  publisher={Wiley Online Library}
}

@article{chen2020predicting,
  title={Predicting station level demand in a bike-sharing system using recurrent neural networks},
  author={Chen, Po-Chuan and Hsieh, He-Yen and Su, Kuan-Wu and Sigalingging, Xanno Kharis and Chen, Yan-Ru and Leu, Jenq-Shiou},
  journal={IET Intelligent Transport Systems},
  volume={14},
  pages={554--561},
  year={2020},
  publisher={Wiley Online Library}
}

@article{kim2019graph,
  title={Graph convolutional network approach applied to predict hourly bike-sharing demands considering spatial, temporal, and global effects},
  author={Kim, Tae San and Lee, Won Kyung and Sohn, So Young},
  journal={PloS One},
  volume={14},
  pages={e0220782},
  year={2019},
  publisher={Public Library of Science San Francisco, CA USA}
}

@inproceedings{yoon2012cityride,
  title={Cityride: a predictive bike sharing journey advisor},
  author={Yoon, Ji Won and Pinelli, Fabio and Calabrese, Francesco},
  booktitle={2012 IEEE 13th International Conference on Mobile Data Management},
  pages={306--311},
  year={2012},
  organization={IEEE}
}

@article{almannaa2020dynamic,
  title={Dynamic linear models to predict bike availability in a bike sharing system},
  author={Almannaa, Mohammed H and Elhenawy, Mohammed and Rakha, Hesham A},
  journal={International Journal of Sustainable Transportation},
  volume={14},
  pages={232--242},
  year={2020},
  publisher={Taylor \& Francis}
}

@article{cantelmo2020low,
  title={Low-dimensional model for bike-sharing demand forecasting that explicitly accounts for weather data},
  author={Cantelmo, Guido and Kucharski, Rafa{\l} and Antoniou, Constantinos},
  journal={Transportation Research Record},
  volume={2674},
  pages={132--144},
  year={2020},
  publisher={SAGE Publications Sage CA: Los Angeles, CA}
}

@article{jiang2022bike,
  title={Bike sharing usage prediction with deep learning: a survey},
  author={Jiang, Weiwei},
  journal={Neural Computing and Applications},
  volume={34},
  pages={15369--15385},
  year={2022},
  publisher={Springer}
}

@inproceedings{mrazovic2018deep,
  title={A deep learning approach for estimating inventory rebalancing demand in bicycle sharing systems},
  author={Mrazovic, Petar and Larriba-Pey, Josep Luis and Matskin, Mihhail},
  booktitle={2018 IEEE 42nd Annual Computer Software and Applications Conference (COMPSAC)},
  volume={2},
  pages={110--115},
  year={2018},
  organization={IEEE}
}

@article{zhou2018markov,
  title={A Markov chain based demand prediction model for stations in bike sharing systems},
  author={Zhou, Yajun and Wang, Lilei and Zhong, Rong and Tan, Yulong and others},
  journal={Mathematical Problems in Engineering},
  volume={2018},
  pages={8028714},
  year={2018},
  publisher={Hindawi}
}

@inproceedings{yang2016mobility,
  title={Mobility modeling and prediction in bike-sharing systems},
  author={Yang, Zidong and Hu, Ji and Shu, Yuanchao and Cheng, Peng and Chen, Jiming and Moscibroda, Thomas},
  booktitle={Proceedings of the 14th Annual International Conference on Mobile Systems, Applications, and Services},
  pages={165--178},
  year={2016}
}

@article{rixey2013station,
  title={Station-level forecasting of bikesharing ridership: Station network effects in three {US} systems},
  author={Rixey, R Alexander},
  journal={Transportation Research Record},
  volume={2387},
  pages={46--55},
  year={2013},
  publisher={SAGE Publications Sage CA: Los Angeles, CA}
}

@article{boufidis2020development,
  title={Development of a station-level demand prediction and visualization tool to support bike-sharing systems’ operators},
  author={Boufidis, Neofytos and Nikiforiadis, Andreas and Chrysostomou, Katerina and Aifadopoulou, Georgia},
  journal={Transportation Research Procedia},
  volume={47},
  pages={51--58},
  year={2020},
  publisher={Elsevier}
}

@inproceedings{ranaiefar2016bike,
  title={Bike sharing ridership forecast using structural equation modeling},
  author={Ranaiefar, Fatemeh and Rixey, R Alexander},
  booktitle={Transportation Research Board 95th Annual Meeting},
  number={16-6573},
  year={2016}
}

@article{liu2019multi,
  title={Multi features and multi-time steps {LSTM} based methodology for bike sharing availability prediction},
  author={Liu, Xu and Gherbi, Abdelouahed and Li, Wubin and Cheriet, Mohamed},
  journal={Procedia Computer Science},
  volume={155},
  pages={394--401},
  year={2019},
  publisher={Elsevier}
}

@article{pan2019predicting,
  title={Predicting bike sharing demand using recurrent neural networks},
  author={Pan, Yan and Zheng, Ray Chen and Zhang, Jiaxi and Yao, Xin},
  journal={Procedia Computer Science},
  volume={147},
  pages={562--566},
  year={2019},
  publisher={Elsevier}
}

@article{kaltenbrunner2010urban,
  title={Urban cycles and mobility patterns: Exploring and predicting trends in a bicycle-based public transport system},
  author={Kaltenbrunner, Andreas and Meza, Rodrigo and Grivolla, Jens and Codina, Joan and Banchs, Rafael},
  journal={Pervasive and Mobile Computing},
  volume={6},
  pages={455--466},
  year={2010},
  publisher={Elsevier}
}

@article{chiariotti2020bike,
  title={A bike-sharing optimization framework combining dynamic rebalancing and user incentives},
  author={Chiariotti, Federico and Pielli, Chiara and Zanella, Andrea and Zorzi, Michele},
  journal={ACM Transactions on Autonomous and Adaptive Systems (TAAS)},
  volume={14},
  pages={1--30},
  year={2020},
  publisher={ACM New York, NY, USA}
}

@article{caggiani2013dynamic,
  title={A dynamic simulation based model for optimal fleet repositioning in bike-sharing systems},
  author={Caggiani, Leonardo and Ottomanelli, Michele},
  journal={Procedia-Social and Behavioral Sciences},
  volume={87},
  pages={203--210},
  year={2013},
  publisher={Elsevier}
}

@article{angelelli2022simulation,
  title={A simulation framework for a station-based bike-sharing system},
  author={Angelelli, Enrico and Chiari, Massimo and Mor, Andrea and Speranza, M Grazia},
  journal={Computers \& Industrial Engineering},
  volume={171},
  pages={108489},
  year={2022},
  publisher={Elsevier}
}

@article{sohrabi2021dynamic,
  title={Dynamic bike sharing traffic prediction using spatiotemporal pattern detection},
  author={Sohrabi, Soheil and Ermagun, Alireza},
  journal={Transportation Research Part D: Transport and Environment},
  volume={90},
  pages={102647},
  year={2021},
  publisher={Elsevier}
}

@article{gammelli2022predictive,
  title={Predictive and prescriptive performance of bike-sharing demand forecasts for inventory management},
  author={Gammelli, Daniele and Wang, Yihua and Prak, Dennis and Rodrigues, Filipe and Minner, Stefan and Pereira, Francisco Camara},
  journal={Transportation Research Part C: Emerging Technologies},
  volume={138},
  pages={103571},
  year={2022},
  publisher={Elsevier}
}

@article{shui2020review,
  title={A review of bicycle-sharing service planning problems},
  author={Shui, CS and Szeto, WY},
  journal={Transportation Research Part C: Emerging Technologies},
  volume={117},
  pages={102648},
  year={2020},
  publisher={Elsevier}
}
\end{document}